\documentclass[graybox]{svmult}
\usepackage{mathptmx}       
\usepackage{helvet}         
\usepackage{courier}        
\usepackage{type1cm}        
\usepackage{makeidx}         
\usepackage{graphicx}        
\usepackage{multicol}        
\usepackage[bottom]{footmisc}

\usepackage{hyperref}
\usepackage{amssymb, mathrsfs,amsmath}
\usepackage{latexsym,array}
\usepackage{amsfonts}
\usepackage{shadow}

\newtheorem{Pa}{Paper}[section]
\newtheorem{Tm}[Pa]{{\bf Theorem}}
\newtheorem{La}[Pa]{{\bf Lemma}}
\newtheorem{Dn}[Pa]{{\bf Definition}}

\newtheorem{Pn}[Pa]{{\bf Proposition}}

\newtheorem{Ex}[Pa]{{\bf Example}}

\newtheorem{Hyp}[Pa]{{\bf Hypothesis}}
\def\e{\epsilon_N}
\def\C{\mathbb C}
\date{}
\begin{document}
\title*{
Extending wavelet filters. Infinite dimensions, the non-rational case,
and indefinite-inner product spaces}
\titlerunning{Extending
wavelet filters}
\author{Daniel Alpay, Palle Jorgensen and Izchak Lewkowicz}
\institute{Daniel Alpay \at Department of Mathematics, Ben Gurion
University of the Negev,  P.O.B. 653, Be'er Sheva 84105,
Israel\email{dany@math.bgu.ac.il} \and Palle Jorgensen\at
Department of Mathematics, 14 MLH, The University of Iowa Iowa
City, IA 52242-1419 USA \email{jorgen@math.uiowa.edu} \and Izchak
Lewkowicz \at Department of Electrical Engineering, Ben Gurion
University of the Negev, P.O.B. 653, Be'er Sheva 84105, Israel
 \email{izchak@ee.bgu.ac.il}}

\maketitle

\abstract{
%
%
In this paper we are discussing various aspects of wavelet
filters. While there are earlier studies of these filters as
matrix valued functions in wavelets, in signal processing, and in
systems, we here expand the framework. Motivated by applications,
and by bringing to bear tools from reproducing kernel theory, we
point out the role of non-positive definite Hermitian inner
products (negative squares), for example Krein spaces, in the
study of stability questions. We focus on the non-rational case,
and establish new connections with the theory of generalized
Schur functions and their associated reproducing kernel
Pontryagin spaces, and the Cuntz relations.
%
%
}

\mbox{}
\section{Introduction}
\setcounter{equation}{0} \label{sec:1}
%
Roughly speaking, systems whose inputs and outputs may be viewed
as {\sl signals} are called {\sl filters}. Mathematically,
filters are often presented as operator valued functions of a
complex variable. In applications, filters are used in areas as
$(i)$ prediction, $(ii)$ signal processing, $(iii)$ systems
theory and $(iv)$ Lax-Phillips scattering theory \cite{LP89}.
There, one is faced with spectral theoretic questions which can
be formulated and answered with the use of a suitable choice of
an operator valued function defined on a domain in complex plane;
in the case of scattering theory, the scattering operator and the
scattering matrix; in the other areas, the names used include
polyphase matrix, see e.g., \cite{Ke04, Kro09}. We also mention
that more recently, filters are used in $(iv)$ multiresolution
analysis in wavelets. We follow standard conventions regarding
time-frequency duality, i.e., the correspondence between discrete
time on one side and a complex frequency variable on the other.
In the simplest cases, one passes from a time series to a
generating function of a complex variable. These frequency
response functions fall in various specific classes of functions
of a complex variable; the particular function spaces in turn are
dictated by applications. Again, motivated by applications, in
our present study, we adopt a wider context for both sides of the
duality divide. On the frequency side, we work with operator
valued functions. This framework is relevant to a host of
applications, and we believe of independent interest in operator
theory. From the literature, we mention \cite{MuSo08}, \cite{CMS}
(see also \cite{2010arXiv1003.2806C}),
and the papers referenced below.\\

We here consider the set of $\mathbb C^{N\times N}$-valued
functions meromorphic in the open unit disk $\mathbb D$ \footnote{Classically, 
in
the engineering literature, the functions are analytic, or more generally 
meromorphic, outside the closed unit disk. The map $z\mapsto 1/z$
relates the two settings.} and
define two subsets of it: We shall denote by $\mathscr C_N$ the
family satisfying the symmetry
\begin{equation}
W(\e z)=W(z)P_N, \label{eq:sym11}
\end{equation}
where $\e=e^{\frac{2\pi i}{N}}$ and  $P_N$ denotes the
permutation matrix
\begin{equation}
\label{Cardinal_Lemoine_Ligne_10}
P_N=\begin{pmatrix} 0_{1\times (N-1)}&1\\
I_{N-1}& 0_{(N-1)\times 1}\end{pmatrix}.
\end{equation}
We shall also denote by $\mathscr U^{I_N}$ the set of
$\mathbb C^{N\times N}$-valued
functions which take unitary values \footnote{For rational
functions, the term {\sl para-unitary} is also used in the
engineering literature.} on the unit circle $\mathbb T$. \vskip
0.2cm

Classically ~{\em wavelet filters}, denoted by $\mathscr W_N$,
are characterized by rational functions satisfying both
symmetries, i.e.
\begin{equation}\label{eq:WN}
\mathscr W_N=\mathscr U^{I_N}\cap\mathscr C_N.
\end{equation}
In a previous paper, see \cite{ajl1}, we have provided an
easy-to-compute characterization of $\mathscr W_N$ as both a set
of rational functions, and in terms of state space realization.
\vskip 0.2cm

The aim of this work is to explore the possibility of extending
the notion of wavelet filters, described in \eqref{eq:WN}. The
functions considered still satisfy the symmetry in
\eqref{eq:sym11}, but:
\begin{itemize}
\item{}The functions are not necessarily rational
or finite dimensional.

\item{}The functions are not necessarily unitary on
the unit circle $\mathbb T$.

\item{}The functions are meromorphic (rather than analytic)
in $\mathbb D$.
\end{itemize}
\vskip 0.2cm

To explain our strategy, first recall the following: If $W$ is a
$\mathbb C^{N\times N}$-valued function which is rational and
takes unitary values on the unit circle, the kernel
\[
K_W(z,w)=\frac{I_N-W(z)W(w)^*}{1-zw^*}
\]
is positive definite in the open unit disk $\mathbb D$ if $W$ has
no poles there, or more generally has a finite number of negative
squares in $\mathbb D$. See Definition \ref{def:nsq} below for the
latter. In our approach, unitarity on the unit circle is replaced
by the requirement that $W$ is a generalized Schur function, in
the sense that $W$ is meromorphic in $\mathbb D$ and the
associated kernel $K_W(z,w)$ has a finite number of negative
squares there. {\sl This family includes in particular the case of
matrix-valued rational functions which take contractive values on
the unit circle}. We will also consider the case where the values
on the unit circle are, when defined, contractive with respect to
indefinite metrics. These kernels are of the form
\begin{equation}
\label{nonsquare} \frac{J_2-W(z)J_1W(w)^*}{1-zw^*}
\end{equation}
when $W$ is $\mathbb C^{p_2\times p_1}$-valued and analytic in a
neighborhood of the origin, and where $J_1$ and $J_2$ are
signature matrices, respectively in ${\mathbb C}^{p_1\times p_1}$
and  ${\mathbb C}^{p_2\times p_2}$, which have the same number of
strictly negative eigenvalues:
\begin{equation}
\label{J1J2} \nu_-(J_1)=\nu_-(J_2),
\end{equation}
and such that the kernel $K_W$ has a finite number of negative
squares. In \cite{ajl1} we studied the realization of wavelet
filters in the $\mathbb C^{N\times M}$-valued (with $M\ge N$)
rational case. The above approach allows us to extend these
results to the case where the filter is not necessarily rational
and $M$ may be smaller than $N$. Furthermore, the conditions in
\cite{ajl1} of the function being analytic in the open unit disk,
and taking coisometric values on the unit circle, are both
relaxed (in particular, in the previous
case, in \eqref{J1J2}, we had $J_1=I_M$ and $J_2=I_N$).\\

The paper is organized as follows. Since we address different
audiences,  Sections 2,3 and 4 are of a review nature. In Section
\ref{sec:2}, we give background on the use of filters in
mathematics. We note that the more traditional framework in the
literature has so far been unnecessarily restricted by two kinds
of technical assumptions: $(i)$ restricting to rational operator
valued functions, and $(ii)$ restricting the range of the operator
valued functions considered. In Section \ref{sec:3} we address
indefinite inner product spaces, and survey the theory of
Pontryagin and Krein spaces. This overview allows us in Section
\ref{sec:4} to describe a setting that expands both the above
mentioned restrictions in (i) and (ii), namely the theory of
generalized Schur functions. Our results in Sections \ref{sec:5}
and \ref{sec:6} (Theorems \ref{th:cuntz12}, \ref{tm:gleason}, and
\ref{tm:dec_rkhs}) deal with representations. We use these
results in obtaining classifications, and decomposition theorems.
In Section \ref{sec:7}, we employ these theorems in the framework
of wavelets.

\section{Some background}
\label{sec:2}

\subsection{Cuntz relations}
The Cuntz relations were realized by J. Cuntz in \cite{Cun77} as
generators of a simple purely infinite $C^*$-algebra. Since then,
they found many applications, and the related literature about
Cuntz relations has flourished. Since Cuntz's paper \cite{Cun77},
the study of their representations has mushroomed, and now makes
up a big literature, see for example \cite{BrJo97, BrJo02a,
BrJo02b,DMP08, BJMP05, Jor03, Jor06c}, and some of their
applications \cite{Jor06a, Jor06b, Jor08, JoSo09, DuJo06a}, for
example to fractals \cite{DuJo06b}.

In the initial framework, one is given a finite set $S_1,\ldots,
S_N$ of isometries with orthogonal ranges adding up to the whole
Hilbert space. Their representations play a role in a variety of
applications, for example wavelets, and more generally
multi-scale phenomena. The study of what are called non-type $I$
$C^*$-algebras was initiated in the pioneering work of Glimm
\cite{Gli61a, Gli61b} and Dixmier \cite{MR0458185}. This in turn
was motivated by use of direct integrals in representation
theory, both in the context of groups and $C^*$-algebras. Direct
integrals of representations are done practically with the use of
Borel cross sections. Glimm proved that there are purely infinite
$C^*$-algebras which do not admit Borel cross sections as a
parameter space for the set of equivalence classes of irreducible
representations; the Cuntz algebra(s) $O_N$  is the best known
examples, \cite{Cun77}. Nonetheless, it was proved in
\cite{BrJo02a} that there are families of equivalence classes of
representations of $O_N$ indexed by wavelet filters, the latter
in turn being indexed by
infinite-dimensional groups.\\
One illustration of the need for expanding the framework of $O_N$
from Hilbert space to the case of Krein spaces is illustrated by
applications to scattering theory for the automorphic wave
equation \cite{LP76}. The initial study was restricted to the
case when the operators $S_i$  act on Hilbert space, and when
they act isometrically. However, since then, there has been a
need for generalizing the Cuntz relations. It was noted in
\cite{BrJo97} that the isometric case adapts well to the
restricted framework of orthogonal wavelet families \cite{Dau92}.
Nonetheless, applications to engineering dictate much wider
families, such as wavelet frames.\\

  {\sl In this work we extend what is known in
the literature in a number of different directions, including to
the case of Pontryagin spaces. We obtain Cuntz relations for
isometries between certain reproducing kernel Pontryagin spaces
of analytic functions.}

\subsection{Wavelet filters}

In electrical engineering terminology, systems whose inputs and
outputs may be viewed as {\sl signals} are called {\sl filters}.
By filter, we here mean functions $W(z)$ defined on the disk in
the complex plane and taking operator
values, i.e., linear operators mapping between suitable spaces.\\

While filters (in the sense of systems and signal processing)
have already been used with success in analysis of wavelets, so
far some powerful tools from systems theory have not yet been
brought to bear on wavelet filters. The traditional restriction
placed on these functions $W(z)$ is that they are rational, and
take values in the unitary group when $z$ is restricted to have
modulus $1$. In models from systems theory, the complex variable
$z$ plays the role of complex frequency. A reason for the recent
success of wavelet algorithms is a coming together of tools from
engineering and harmonic analysis. While wavelets now enter into
a multitude of applications from analysis and probability, it was
the incorporation of ideas from signal processing that offered
new and easy-to-use algorithms, and hence wavelets are now used
in both discrete problems, as well as in harmonic analysis
decompositions. It is our purpose to use tools from systems
theory in wavelet problems and also show how ideas from wavelet
decompositions shed light on factorizations used by engineers.
Each of the various wavelet families demands a separate class of
filters, for the case of compactly supported biorthogonal
wavelets, see for example Resnikoff, Tian, Wells \cite{MR1858875}
and Sebert and Zou \cite{2011arXiv1101.3793S}. By now there is a
substantial literature on the use of filters in wavelets (see
e.g., \cite{BrJo02a, Dau92, Jor03, Jor06c}). For filters in
wavelets, there are two pioneering papers \cite{Law91a, Law91b},
and the book \cite{Mal09}.\\

In a previous work \cite{ajl1} we characterized all rational
wavelet filters attaining unitary values on the unit circle. It
turned out that this family is quite small ( and in particular
the subset of Finite Impulse Response filters, commonly used in
engineering).\\

{\sl Thus, we here remove both restrictions on the filters, i.e.,
rational and unitary, and consider $W(z)$ which are generalized
Schur functions, and  use reproducing kernel Pontryagin spaces
associated with $W$. See \cite{adrs} for background.}\\

We hope that this message will be useful to practitioners in
their use of these rigorous mathematics tools.

\section{Pontryagin spaces and Krein spaces}
\label{sec:3} For a number of problems
in the study of signals and filters (for example stability
considerations), it is necessary to work with Hermitian inner
products that are not positive definite. This view changes the
Hermitian quadratic forms, allowing for negative squares, as well
as the associated linear spaces. But more importantly, this wider
setting also necessitates changes in the analysis, for example in
the meaning of the notion of the adjoint operator, as well as the
reproducing kernels. There are a number of subtle analytic points
involved, as well as a new operator theory. We turn to these
details below.

\subsection{Krein spaces}
  A {\sl Krein space} is a pair $(V, [\cdot,\cdot])$, where $V$ is a
linear vector space on $\mathbb C$ endowed with an Hermitian form
$[ \cdot,\cdot]$, and with the following properties: $V$ can be
written as $V=V_++V_-$, where:
\begin{enumerate}
\item
$V_+$ endowed with the Hermitian  form $[ \cdot,\cdot]$ is a
Hilbert space.

\item
$V_-$ endowed with the Hermitian form $-[ \cdot,\cdot]$ is a
Hilbert space.

\item It holds that $V_+\cap V_-=\left\{0\right\}$.

\item
For all $v_\pm\in V_\pm$,
\[
[ v_+,v_-]=0.
\]
\end{enumerate}
The representation $V=V_++V_-$ is called a {\sl fundamental
decomposition}, and is highly non unique as soon as ${\rm
dim}~V_->0$. Given such a decomposition, the map
\[
\sigma(v_++v_-)=v_+-v_-
\]
is called a {\sl fundamental symmetry}. Note that the space $V$ endowed
with the Hermitian form (where $w=w_++w_-$ is also an element of
$V$, with $w_\pm\in V_\pm$)
\[
\langle v,w\rangle=[v,\sigma w]=[v_+,w_+]-[v_-,w_-]
\]
is a Hilbert space. These norms are called {\sl natural norms},
and they are all equivalent. The Hilbert space topologies
associated to any two such decompositions are equivalent, and $V$
is endowed with any of them; see \cite[p. 102]{bognar}. When
$V_-$ is finite dimensional, $V$ is called a Pontryagin space and
the dimension of $V_-$ is called the {\sl negative index} (or the
{\sl index} for short) of the Pontryagin space. We refer to the
books \cite{azih}, \cite{bognar}, \cite{ikl}, \cite{adrs} for more
information on Krein and Pontryagin spaces. Note that in
\cite{ikl} it is the space $V_+$ rather than $V_-$ which is
assumed finite dimensional in the definition of a Pontryagin
space. Surveys may be found in for instance in \cite{dj-field},
\cite{MR92m:47068}, \cite{ad3}. It is interesting to note that
Laurent Schwartz introduced independently the notion of Krein and
Pontryagin spaces (he used the terminology Hermitian spaces for
Krein and Pontryagin spaces) in his paper \cite{schwartz}. 
For applications of Krein spaces to the study of boundary conditions 
for hyperbolic PDE, including wave equations, and exterior domains, 
see for example \cite{CP68, LaPh70, LaPh71, Phi61}. We now
give two examples, which will be important in the sequel.

\begin{Ex} Let $J\in\mathbb C^{p\times p}$ be an Hermitian
involution, i.e.
\[
J=J^{-1}=J^*.
\]
Such a matrix is called a signature matrix. We denote by $\mathbb
C_J$ the space $\mathbb C^p$ endowed with the associated
indefinite inner product
\[
[x,y]_J=y^*Jx,\quad x,y\in\mathbb C^p.
\]
It is a finite dimensional Pontryagin space. \label{paris_texas}
\end{Ex}

\begin{Ex}
Let $J$ be a signature matrix. We consider the space $\mathbf
H_{2}(\mathbb D)^p$ of functions analytic in $\mathbb D$ and with
values in $\mathbb C^p$:
\[
f(z)=\sum_{n=0}^\infty a_n z^n,\quad a_n\in\mathbb C^p,
\]
such that
\[
\sum_{n=0}^\infty a_n^*a_n<\infty.
\]
Then, $\mathbf H_2(\mathbb D)^p$ endowed with the Hermitian form
\[
[ f,g]_J=\sum_{n=0}^\infty b_n^*Ja_n\quad (with\,\,
g(z)=\sum_{n=0}^\infty b_nz^n)
\]
is a Krein space, which we denote by $\mathbf H_{2,J}(\mathbb D)$.
\end{Ex}

In the above example, if $p=1$ and $J=1$ (as opposed to $J=-1$)
the space $\mathbf H_{2,J}(\mathbb D)$ is equal to the classical
Hardy space $\mathbf H_2(\mathbb D)$ of the unit disk.

\subsection{Operators in Krein and Pontryagin spaces}
When one considers a bounded operator $A$ between two Krein spaces
$(\mathcal K_1,[\cdot ,\cdot]_1)$ and $(\mathcal K_2,[\cdot
,\cdot]_2)$ (in this paper, it will be most of the time between
two Pontryagin spaces) the adjoint can be computed in two
different ways, with respect to the Hilbert spaces inner
products, (and then we use the notation $A^*$) and with respect
to the Krein spaces inner products (and then we use the notation
$A^{[*]}$). More precisely, if $\sigma_1$ and $\sigma_2$ are
fundamental symmetries in $\mathcal K_1$ and $\mathcal K_2$ which
define the Hilbert spaces inner products
\[
\langle f_1,g_1\rangle_1=[\sigma_1f_1,g_1]_1\quad{\rm and}\quad
\langle f_2,g_2\rangle_2=[\sigma_2f_2,g_2]_2,
\]
(with $f_1,g_1\in\mathcal K_1$ and $f_2,g_2\in\mathcal K_2$), we
have for $f_1\in\mathcal K_1$ and $f_2\in\mathcal K_2$
\[
\begin{split}
[Af_1,f_2]_2&=\langle \sigma_2 Af_1,f_2\rangle_2\\
&=\langle f_1,A^*\sigma_2f_2\rangle_1\\
&=[f_1,A^{[*]}f_2]_1,
\end{split}
\]
with
\begin{equation}
A^{[*]}=\sigma_1A^*\sigma_2.
\label{paris}
\end{equation}
In the case of $\mathbb C_J$ (see Example \ref{paris_texas}) we
have
\begin{equation}
A^{[*]}=JA^*J.
\end{equation}

The operator $A$ from $\mathcal D(A)\subset\mathcal K_1$, where
$(\mathcal K_1,[\cdot ,\cdot]_1)$ is a Krein space, into the
Krein space $(\mathcal K_2,[\cdot ,\cdot]_2)$ is a contraction if
\[
[Ak_1,Ak_1]_2\le [k_1,k_1]_1,\quad\forall k_1\in\mathcal D(A).
\]
A densely defined contraction, or even isometry, operator $A$
between Krein spaces need not be continuous, let alone have a
continuous extension. See for instance \cite[Theorem
1.1.7]{MR92m:47068}. In the case of Pontryagin spaces with same
negative index, $A$ has a continuous extension to all of $\mathcal
K_1$, see \cite[Theorem 1.4.1, p. 27]{adrs}, and Theorem
\ref{tm:schmulyan} below. Even when it is continuous and has a
well-defined adjoint, this adjoint need not be a contraction. The
operator is called a bicontraction if both it and its adjoint are
contractions. When the Krein spaces are Pontryagin spaces with
same negative index, a contraction is automatically continuous and
its adjoint is also a contraction. An important notion in the
theory of Pontryagin spaces is that of {\sl relation}. Given two
Pontryagin spaces $\mathcal P_1$ and $\mathcal P_2$, a relation
is a linear subspace of $\mathcal P_1\times\mathcal P_2$. For
instance the graph of an operator is a relation. The domain of
the relation $\mathscr R$ is the set of $f\in\mathcal P_1$ such
that there is a $g\in\mathcal P_2$ for which $(f,g)\in\mathscr
R$. A relation $\mathscr R$ is called {\sl contractive} if,
\[
[g,g]_2\le [f,f]_1\quad\forall (f,g)\in\mathscr R.
\]
A key result is the following theorem of Shmulyan (see
\cite[Theorem 1.4.1, p. 27]{adrs}).
\begin{Tm}
\label{tm:schmulyan} A densely defined contractive relation
between Pontryagin spaces with same negative index extends to the
graph of a uniquely defined contraction operator from $\mathcal
P_1$ into $\mathcal P_2$.
\end{Tm}
\subsection{Kernels}
Recall that a (say, matrix-valued) function $K(z,w)$ of two
variables, defined for $z$ and $w$ in a set $\Omega$ is called a
positive definite kernel if it is Hermitian: $K(z,w)^*=K(w,z)$
for all $z,w\in\Omega$, and if for every choice of $M\in\mathbb
N$ and $w_1,\ldots, w_M\in\Omega$ the $M\times M$ Hermitian block
matrix with $(\ell,j)$ block entry $K(w_\ell,w_j)$ is non
negative. For instance, if $b$ is a finite Blaschke product,
\[
b(z)=\prod_{n=1}^m\frac{z-a_n}{1-za_n^*}
\]
for some $a_1,\ldots,a_m$ in the open unit disk, the kernel
\[
k_b(z,w)=\frac{1-b(z)b(w)^*}{1-zw^*}
\]
is positive definite, as can be seen from the formula
\[
k_b(z,w)=\langle k_b(\cdot, w),k_b(\cdot,z)\rangle_{\mathbf
H_2(\mathbb D)}.
\]
When $b$ is replaced with a function $s$ analytic and contractive
in the open unit disk, the corresponding kernel
$k_s(z,w)=\frac{1-s(z)s(w)^*}{1-zw^*}$ is still positive definite
in $\mathbb D$, see \cite{dbr1}, \cite{dbr2}. This follows, for
instance, from the fact that the operator of multiplication by
$s$ is a contraction from ${\mathbf H_2(\mathbb D)}$ into itself.
In the special case of a finite Blasckhe product (or more
generally, of an inner function), this multiplication operator is
an isometry. This makes the underlying computations much easier.
More generally, the kernels which appear in the following section
can be seen as far
reaching generalizations of the kernels $k_b(z,w)$.\\

The notion of positive definite kernel has been extended by Krein
as follows:

\begin{Dn}
\label{def:nsq} Let $\kappa\in\mathbb N_0$. A (say, matrix-valued)
function $K(z,w)$ defined on a set $\Omega$ has $\kappa$ negative
squares if it is Hermitian, and if for every choice of
$M\in\mathbb N$ and $w_1,\ldots, w_M\in\Omega$ the $M\times M$
Hermitian block matrix with $(\ell,j)$ block entry
$K(w_\ell,w_j)$ has at most $\kappa$ strictly negative
eigenvalues, and exactly $\kappa$ strictly negative eigenvalues
for some choice of $M,w_1,\ldots, w_M$. When $\kappa=0$, the
function is positive definite.
\end{Dn}

The one-to-one correspondence between positive definite kernels
and reproducing kernel Hilbert spaces was first extended to the
indefinite case by L. Schwartz; see \cite{schwartz}: {\sl There
is a one-to-one correspondence between reproducing kernel
Pontryagin spaces and kernels with a finite number of negative
squares}. For completeness, we mention that such a result fails
if the number of negative squares is not finite. {\sl A necessary
and sufficient condition for a function to be the reproducing
kernel of a Krein space is that this function is the difference
of two positive kernels, but the associated Krein space need not
be unique}. Here too we refer to Schwartz \cite{schwartz}, and
also to the paper \cite{a2}. Realization of operator-valued
analytic functions (without assumptions on an associated kernel,
but with some symmetry hypothesis) have also been considered. See for
instance \cite{MR903068}.
The $\mathbb C^{p\times p}$-valued function $K(z,w)$ defined for
$z,w$ in an open set $\Omega$ of the complex plane will be called
an {\sl analytic kernel} if it is Hermitian and if it is analytic in
$z$ and $w^*$. If it has moreover a finite number of negative
squares, the elements of the associated reproducing kernel
Pontryagin space are analytic in $\Omega$.
See \cite[Theorem 1.1.2, p. 7]{adrs}.\\

There are two important classes of operators between reproducing
kernel spaces, namely multiplication and composition operators.
We conclude this section with three results on these operators.

\begin{Tm}
Let $(\mathcal K_1,[\cdot ,\cdot]_1)$ and $(\mathcal K_2,[\cdot
,\cdot]_2)$ be two reproducing kernel Krein spaces of
vector-valued functions, defined in $\Omega$, and with
reproducing kernels $K_1(z,w)$ and $K_2(z,w)$, respectively
$\mathbb C^{p_1\times p_1}$- and $\mathbb C^{p_2\times
p_2}$-valued. Let $m$ be a $\mathbb C^{p_2\times p_1}$-valued
function and let $\varphi$ be a map from $\Omega$ into itself.
Assume that the map
\begin{equation}
\label{tm} (T_{m,\varphi} f)(z)=m(z)f(\varphi(z))
\end{equation}
defines a bounded operator from $(\mathcal K_1,[\cdot ,\cdot]_1)$
into $(\mathcal K_2,[\cdot ,\cdot]_2)$. Then, for every
$z,w\in\Omega$, and $\xi_2\in\mathbb C^{p_2}$,
\begin{equation}
\label{placedelodeon} \left(T_{m,\varphi}^{[*]}K_2(\cdot,
w)\xi_2\right)(z)=K_1(z,\varphi(w))m(w)^*\xi_2.
\end{equation}
\end{Tm}

{\bf Proof:} Let $z,w\in\Omega$, $\xi_2\in\mathbb C^{p_2}$ and
$\xi_1\in\mathbb C^{p_1}$. We have
\[
\begin{split}
\xi_1^*\left(T_{m,\varphi}^{[*]}K_2(\cdot, w)\xi_2\right)(z)&=[
T_{m,\varphi}^{[*]}K_2(\cdot, w)\xi_2\, ,\, K_1(\cdot, z)\xi_1]_1\\
&=[
K_2(\cdot, w)\xi_2\, ,\, T_{m,\varphi}( K_1(\cdot, z)\xi_1)]_2\\
&=[K_2(\cdot, w)\xi_2\, ,\, m(\cdot)K_1(\varphi(\cdot), z)\xi_1]_2\\
&=[ m(\cdot)K_1(\varphi(\cdot), z)\xi_1 \, ,\, K_2(\cdot, w)\xi_2
]_2^*\\
&=\left(\xi_2^*m(w)K_1(\varphi(w)\, ,\, z)\xi_1\right)^*\\
&=\xi_1^*K_1(z,\varphi(w))m(w)^*\xi_2.
\end{split}
\]
\mbox{}\qed\mbox{}\\

As a corollary we have the following result:
\begin{Tm}
Assume in the preceding theorem that $\mathcal K_1$ and $\mathcal
K_2$ are Pontryagin spaces with same negative index. Then,
$T_{m,\varphi}$ is a contraction if and only if the kernel
\begin{equation}
\label{montparnasse} K_2(z,w)-m(z)K_1(\varphi(z),\varphi(w))m(w)^*
\end{equation}
is positive definite in $\Omega$. \label{th:contraction}
\end{Tm}

{\bf Proof:} Assume that $T$ is a contraction. Then, its adjoint
is also a contraction since the Pontryagin spaces have the same
negative index. Let $g\in\mathcal K_2$ be of the form
\[
g(z)=\sum_{k=1}^N K_2(z,w_k)\xi_k,
\]
where $N\in\mathbb N$, $w_1,\ldots, w_N\in\omega$ and
$\xi_1,\ldots, \xi_N\in\mathbb C^{p_2}$. By \eqref{placedelodeon}
we have
\[
\begin{split}
\sum_{\ell,k=1}^N\xi_\ell^*m(w_\ell)K_1(\varphi(w_\ell),
\varphi(w_k))m(w_k)^*\xi_k&=\\
&\hspace{-5cm}= [\sum_{k=1}^N K_1(z,\varphi(w_k))m(w_k)^*\xi_k,
\sum_{\ell=1}^N K_1(z,\varphi(w_\ell))m(w_\ell)^*\xi_\ell]_1\\&
\hspace{-5cm}=
[T_{m,\varphi}^{[*]}g,T_{m,\varphi}^{[*]}g]_1\\
&\hspace{-5cm}\le [g,g]_2\\
&\hspace{-5cm}= \sum_{\ell,k=1}^N\xi_\ell^*K_2(w_\ell,w_k)\xi_k,
\end{split}
\]
and hence the kernel \eqref{montparnasse} is positive definite.
Conversely, assume that the kernel \eqref{montparnasse} is
positive definite. Then the linear span of the pairs of functions
\[
(K_2(\cdot,w)\xi\, ,\, K_1(\cdot,\varphi(w))m(w)^*\xi), \quad
w\in\Omega,\,\, \xi\in\mathbb C^{p_2},
\]
defines a linear densely defined contractive relation in $\mathcal K_1\times
\mathcal K_2$. By Shmulyan's theorem (see Theorem
\ref{tm:schmulyan}), this relation has an everywhere defined
extension which is the graph of a bounded contraction: There is a
unique contraction $X$ from $\mathcal K_2$ into $\mathcal K_1$
such that
\[
X(K_2(\cdot,w)\xi)=K_1(\cdot,\varphi(w))m(w)^*\xi,\quad
w\in\Omega,\,\, \xi\in\mathbb C^{p_2}.
\]
By \eqref{placedelodeon}, we have $X^{[*]}=T_{m,\varphi}$, and
this
concludes the proof.\mbox{}\qed\mbox{}\\

We will consider in the sequel special cases of this result, in particular
when
\[
m(z)=\begin{pmatrix}1&z&\cdots&z^{N-1}\end{pmatrix},
\]
see Theorem \ref{th:cuntz12}, or more generally when
\[
m(z)=\begin{pmatrix}m_0(z)&m_1(z)&\cdots&m_{N-1}(z)\end{pmatrix},
\]
see Theorem \ref{tm:gleason}. The operator $T_{m,\varphi}$
defined by \eqref{tm} is then a block operator, and its components
satisfy, under appropriate supplementary hypothesis, the Cuntz
relations formally defined in \eqref{C1}-\eqref{C2} below.\\

We conclude this section with a result on composition operators in
reproducing kernel Pontryagin spaces.

\begin{Tm}
\label{tm:varphi} 
Let $K(z,w)$ be a $\mathbb C^{p\times
p}$-valued function which has $\kappa$ negative squares in the set
$\Omega$. The associated reproducing kernel Pontryagin space will
be denoted by $\mathcal P(K)$. Let $\varphi$ be a map from
$\Omega$ into itself, and assume that:
\[
f(\varphi(z))\equiv0\Longrightarrow f\equiv 0
\]
for $f\in\mathcal P(K)$. Then:\\
$(a)$ The function $K_\varphi(z,w)= K(\varphi(z),\varphi(w))$ has
at most $\kappa$ negative squares in $\Omega$ and its associated
reproducing Pontryagin space is the set of functions of the form
$F(z)= f(\varphi(z))$, with $f\in\mathcal P(K)$ and 
Hermitian form
\begin{equation}
\label{hermitian11} [F,G]_{\mathcal P(K_\varphi)}
=[f,g]_{\mathcal
P(K)}.
\end{equation}
$(2)$ The map $f\mapsto f(\varphi)$ is unitary from $\mathcal
P(K)$ into itself if and only if
\begin{equation}
K(z,w)=K(\varphi(z),\varphi(w)),\quad \forall z,w\in\Omega.
\end{equation}
\end{Tm}

{\bf Proof:} Set
\[
\mathcal M_\varphi=\left\{f(\varphi(z)),\,\, f\in\mathcal
P(K)\right\}.
\]
By hypothesis, we have $f(\varphi(z))\equiv 0$ if and only if
$f\equiv 0$, and so the Hermitian form \eqref{hermitian11} is
well defined and induces a Pontryagin structure on $\mathcal
M_\varphi$. Furthermore, with $c\in\mathbb C^p$ and
$F(z)=f(\varphi(z))\in\mathcal M_\varphi$, we have:
\[
\begin{split}
[F(\cdot), K_\varphi(\cdot,w)c]_{\mathcal P(K_\varphi)}
&=[f(\cdot), K(\cdot,\varphi(w))c]_{\mathcal P(K)}\\
&=c^*f(\varphi(z))\\
&=F(w),
\end{split}
\]
and hence the reproducing kernel property is in force. To prove
$(b)$ we use the uniqueness of the kernel for a given reproducing
kernel Pontryagin space.
\mbox{}\qed\mbox{}\\

To fine-tune the previous result, note that for $\varphi(z)=z^N$,
the composition map is an isometry from $\mathbf H_2(\mathbb D)$
into itself, but is not unitary (unless $N=1$). We also note that
the preceding theorem holds also for reproducing kernel Krein
spaces. Indeed, the correspondence between functions which are
difference of positive functions on a given set and reproducing
kernel Krein spaces is not one-to-one, but a given reproducing
kernel Krein space has a unique reproducing kernel.

\section[Generalized Schur functions]
{Generalized Schur functions and associated spaces}
\label{sec:4} In this section we review
the main aspects of the the realization theory of generalized
Schur functions and of their associated reproducing kernel
Pontryagin spaces.

\subsection{Generalized Schur functions}
\label{sec1_1}

In the positive definite case, this theory originates with the
works of de Branges and Rovnyak, see \cite{dbr1,dbr2}. In earlier
work on models involving operators in Hilbert space, and matrix
factorization, de Branges spaces have served as a surprisingly
powerful tool. The theory was developed in the indefinite case in
in a fundamental series of papers by Krein and Langer, see for instance
\cite{kl1,kl2,MR518342,MR614775,kl3}, and, using reproducing kernel methods
in \cite{ad3} and in the book
\cite{adrs}. It was later used in \cite[p. 119]{adrs} and in the
paper \cite{aadl1} to study generalized Schur functions with some
given symmetry. In this paper we use this setting to present non
rational and non unitary wavelet filters. In \cite{adrs} the case
of operator valued functions is studied, but we here consider the
case of $\mathbb C^{p\times p}$-valued functions. We now recall
the definition of a generalized Schur function. A (say $\mathbb
C^{p\times p}$-valued) function $W$ is called a Schur function if
it is analytic and contractive in the open unit disk, or,
equivalently, if the associated kernel
\begin{equation}
\label{place_d_italie} K_W(z,w)=\frac{I_p-W(z)W(w)^*}{1-zw^*}
\end{equation}
is positive definite in a neighborhood of the origin. Then, it
has a unique analytic extension to the open unit disk, and this
extension is such that the kernel $K_W$ is still positive
definite in $\mathbb D$. There are two other kernels associated
to $W$, namely the kernel $K_{\widetilde{W}}(z,w)$
(with $\widetilde{W}(z)\stackrel{\rm def.}{=}
W(z^*)^*$), and the kernel
\[
D_W(z,w)=\begin{pmatrix} K_W(z,w)&\frac{W(z)-W(w^*)}{z-w^*}\\
\frac{\widetilde{W}(z)-\widetilde{W}(w^*)}{z-w^*}
&K_{\widetilde{W}}(z,w)
\end{pmatrix}.
\]
These three kernels are simultaneously positive definite in the
open unit disk. The first is the state space for a unique
coisometric realization of $W$, the second is the state space for
a unique isometric realization of $W$, and the reproducing kernel
Hilbert space with reproducing kernel $D_W$ is the state space for
a unique unitary realization of $W$. In these three
cases, uniqueness is up to an invertible similarity operator.\\

Let $J\in\mathbb C^{p\times p}$ be a signature matrix. We now
consider functions with values in $\mathbb C_J$ defined in Example
\ref{paris_texas}, denoted by $\Theta$ (rather than $W$).
A $\mathbb C^{p\times p}$-valued functions $\Theta$ analytic in a
neighborhood of the origin is called $J$-contractive if the
associated kernel
\begin{equation}
\label{Ktheta}
K_\Theta(z,w)=\frac{J-\Theta(z)J\Theta(w)^*}{1-zw^*}
\end{equation}
is positive definite. It has a unique meromorphic extension to
the open unit disk, and this extension is such that the kernel
$K_\Theta$ is still positive definite in the domain of
analyticity of $\Theta$ in $\mathbb D$. Here too, besides the
kernel $K_\Theta$ we have the kernel
$K_{\widetilde{\Theta}}(z,w)$ and the kernel
\[
D_\Theta(z,w)=\begin{pmatrix}K_\Theta(z,w)&\frac{J\Theta(z)-J\Theta
(w^*)}{z-w^*}\\
\frac{\widetilde{\Theta}(z)J-\widetilde{\Theta}
(w^*)J}{z-w^*}&K_{\widetilde{\Theta}}(z,w)
\end{pmatrix}.
\]

We note that the kernel $K_\Theta$ can be written as
\[
K_\Theta(z,w)=\frac{I_p-\Theta(z)\Theta(w)^{[*]}}{1-zw^*},
\]
where $[*]$ denotes the adjoint in $\mathbb C_J$. This conforms
with the way these kernels and the two other related kernels are
written down
in \cite{adrs}.\\

As we already mentioned, Krein and Langer developed
in \cite{kl1,kl2,MR518342,MR614775,kl3},
the theory of operator-valued functions such that the
corresponding kernels $K_\Theta$ (with a signature operator
rather than a signature matrix) has a finite number of negative
squares in some open subset of the open unit disk. Then, $\Theta$
has a unique meromorphic extension to the open unit disk, and
this extension is such that $K_\Theta$ has the same number of
negative squares in $\Omega(\Theta)$, the domain of analyticity
$\Theta$ in $\mathbb D$.  The three kernels have simultaneously
the same number of negative squares, and as in the positive
definite case, are respectively state spaces for coisometric,
isometric and
unitary realizations of $\Theta$.\\

In the special case $J=I$  (we return to the notation $W$ rather
than $\Theta$ for the function), Krein and Langer proved, see
\cite{kl1}, that $W$ can be written as $W_0B_0^{-1}$, where $W_0$
is analytic and contractive in the open unit disk, and where
$B_0$ is a finite matrix-valued Blaschke product. It follows that
$W$ has a finite number of poles in the open unit disk. In the
rational case, and when $W$ takes unitary values on the unit
circle, $W$ is a quotient of two matrix-valued rational Blaschke
product. Note however that when $J$ has mixed inertia, $W$ may
have an infinite number of poles, even when $\kappa=0$. For
example, take
\[
J=\begin{pmatrix}1&0\\0&-1\end{pmatrix}\quad{\rm and}\quad
W(z)=\begin{pmatrix}1&0\\0&b(z)^{-1}\end{pmatrix},
\]
where $b$ is a convergent Blasckhe product with an infinite
number of zeros. Such examples originate with the work of Potapov
\cite{pootapov}.

\begin{Dn}
We denote by $\mathscr S^{p\times p}_\kappa(\mathbb D)$ the
family of $\mathbb C^{p\times p}$-valued functions $W$
meromorphic in the open unit disk, and such that the kernel $K_W$
(defined by \eqref{place_d_italie}) has $\kappa$ negative squares
in the domain of
analyticity of $W$ in $\mathbb D$.\\
Given a signature matrix $J$, we denote by $\mathscr
S^J_{\kappa}(\mathbb D)$ the family of $\mathbb C^{p\times
p}$-valued functions $\Theta$ meromorphic in the open unit disk,
and such that the kernel $K_\Theta$ (defined by \eqref{Ktheta})
has $\kappa$
negative squares in the domain of analyticity of $\Theta$ in $\mathbb D$.\\
We denote by $\mathcal P(W)$ and $\mathcal P(\Theta)$
respectively the associated reproducing kernel Pontryagin spaces.
\label{def:sjk}
\end{Dn}

Since the kernels $K_W$ and $K_\Theta$ are analytic in $z$ and
$w^*$, the elements of the associated reproducing kernel
Pontryagin spaces are analytic in the domain of definition of $W$
or $\Theta$ respectively. See \cite[Theorem 1.1.3, p. 7]{adrs}.\\

More generally, it is useful to consider non square generalized
Schur functions. We consider $J_1\in\mathbb C^{p_1\times p_1}$
and $J_2\in\mathbb C^{p_1\times p_1}$ two signature matrices, of
possibly different sizes, such that \eqref{J1J2} is in form
denoted by $\nu_-$:
\begin{equation*}
\nu_-(J_1)=\nu_-(J_2).
\end{equation*}
Reproducing kernel Pontryagin spaces with reproducing kernel of
the form \eqref{nonsquare}:
\begin{equation*}
\frac{J_2-\Theta(z)J_1\Theta(w)^*}{1-zw^*}
\end{equation*}
when $\Theta$ is $\mathbb C^{p_2\times p_1}$-valued and analytic
in a neighborhood of the origin, have been characterized in
\cite[Theorem 3.1.2, p. 85]{adrs} (in fact, the result there is
more general and considers operator-valued functions). In the
statement below $R_0$ denotes the backward-shift operator
\[
R_0f(z)=\frac{f(z)-f(0)}{z}.
\]
\begin{Tm}
Let $(\mathcal P, [\cdot, \cdot]_{\mathcal P})$ be a reproducing
kernel Pontryagin space of $\mathbb C^{p_2}$-valued functions. It
has a reproducing kernel of the form \eqref{nonsquare} if and
only if it is invariant under the backward-shift operator $R_0$
and
\[
[R_0f,R_0f]_{\mathcal P}\le [f,f]_{\mathcal
P}-f(0)^*J_2f(0),\quad\forall f\in\mathcal P.
\]
\end{Tm}

An example of such non square $\Theta$ appears in Section
\ref{azxcv} below. See \eqref{azxc}.\\

\subsection{State spaces and realizations}
\label{sub42}
We begin with recalling the following definition. Let $W$  be an
operator-valued function analytic in a neighborhood of the
origin. A realization of $W$ is an expression of the form
\begin{equation}
\label{realP} W(z)=D+zC(I-zA)^{-1}B,
\end{equation}
where $D=W(0)$ and $A,B,C$ are operators between appropriate
spaces. It is an important problem to connect the properties of
$W$ and of the operator matrix
\begin{equation}
\label{Cambronne} M=\begin{pmatrix}A&B\\ C&D\end{pmatrix}.
\end{equation}
When the values of $W$ are linear bounded operators between two
Krein spaces, Azizov proved that a realization exists, and that
$M$ can be chosen unitary. See \cite{az2}, and see
\cite{MR903068} for further discussion and additional references.
When $W$ is a matrix-valued rational function without a pole at
the origin, the spaces may be chosen finite dimensional, when no
special
structure is forced on the operator matrix $M$.\\

In Section \ref{sec1_1}, we have studied the correspondence
between kernels and operator valued Schur functions. Here we then
pass to the realizations of Schur functions. The introduction of
Schur functions offers many advantages, relevant to algorithms
and to computation. Case in point: In the next subsection, we
give explicit formulas for realizations, i.e., for the
computation of the four block operator entries $A$ through $D$
making up admissible realizations of a given Schur function, and
therefore of a kernel.  As we show, there are several such
choices, the coisometric realization (Theorem \ref{tm:coiso}),
and the unitary realization of de Branges and Rovnyak (Theorem
\ref{tm:unitary}), among others. There is in turn a rich
literature on Schur algorithms in various special cases, see for
example \cite{MR2002b:47144} for an overview. In preparation of
Section \ref{minimal!} we need some definitions. Let $\mathcal P$
denote the space where $A$ acts in \eqref{realP}. We say that the
realization is closely inner connected if the span of the
functions
\[
(I-zA)^{-1}B\xi,
\]
where $\xi$ runs through $\mathbb C^p$ (recall that $J\in\mathbb
C^{p\times p}$) and $z$ runs through a neighborhood of the origin,
is dense in $\mathcal P$. With the same choices of $\xi$ and $z$,
it will be called closely outer connected if the span of the
functions
\[
(I-zA^{[*]})^{-1}C^{[*]}\xi
\]
is dense in $\mathcal P$, and connected if the span of the
functions
\[
(I-zA)^{-1}B\xi,\quad{\rm and}\quad (I-wA^{[*]})^{-1}C^{[*]}\eta
\]
is dense in $\mathcal P$ ($\eta$ running through $\mathbb C^p$
and $w$ through the same neighborhood of the origin as $z$). Here
the adjoints are between Pontryagin spaces. We note that the
terminology is different from that of classical system theory. In
the finite dimensional case, what is called here closely inner
connected corresponds to observability, and what is called outer
connected corresponds to controlabilty. The notion of being
closely connected is specific to this domain, and is, in general,
different from minimality.\\

\subsection{Coisometric and unitary realizations}
\label{minimal!} Let $\Theta\in\mathscr S_\kappa^J$ be a
generalized Schur function, assumed analytic in a neighborhood of
the origin. In this section we review how the spaces $\mathcal
P(\Theta)$ and $\mathcal D(\Theta)$ are the state spaces for
coisometric and unitary realizations respectively. For the
following theorems, see \cite[Theorem 2.2.1, p. 49]{adrs} and
\cite[Theorem 2.1.3]{adrs} respectively. In Theorems
\ref{tm:coiso} and \ref{tm:unitary} below the notions of
coisometry and unitarity means that $M$ in \eqref{Cambronne}
is an operator coisometric (resp. unitary) from the Pontryagin
\mbox{$\mathcal P(\Theta)\oplus \mathbb C_J$} into itself (resp.
from $\mathcal D(\Theta)\oplus \mathbb C_J$ into itself).

\begin{Tm}
\label{tm:coiso} Let $J\in\mathbb C^{p\times p}$ be a signature
matrix, and  $\Theta\in\mathscr S_\kappa^J$ be analytic in a
neighborhood of the origin. Then the formulas
\[
\begin{split}
Af(z)&=\frac{f(z)-f(0)}{z},\\
(B\xi)(z)&=\frac{\Theta(z)-\Theta(0)}{z}\xi,\\
Cf&=f(0),\\
D\xi&=\Theta(0)\xi,
\end{split}
\]
with $f\in\mathcal P(\Theta)$ and $\xi\in\mathbb C^p$, define a
closely outer connected realization of $\Theta$ which is
coisometric. This realization is unique up to a continuous and
continuously invertible similarity operator.
\end{Tm}

This coisometric realization was introduced by L. de Branges and
J. Rovnyak in \cite{dbr1} for scalar Schur functions, and
extended to the operator-valued case in \cite{dbr2}. We note that
the coisometric realization is also known as the {\sl backward
shift realization};
see e.g. \cite{MR1393938}.\\

L. de Branges and J. Rovnyak also formulated the unitary
realization below.
\begin{Tm}
Let $J\in\mathbb C^{p\times p}$ be a signature matrix, and
$\Theta\in\mathscr S_\kappa^J$ be analytic in a neighborhood of
the origin. The formulas
\[
\begin{split}
A\begin{pmatrix}f\\g\end{pmatrix}&=\begin{pmatrix}
\dfrac{f(z)-f(0)}{z}\\
zg(z)-\widetilde{\Theta}(z)Jf(0)\end{pmatrix},\\
(B\xi)(z)&=\begin{pmatrix}\dfrac{\Theta(z)-\Theta(0)}{z}\xi\\
(J-
\widetilde{\Theta}(z)J\widetilde{\Theta}(0)^*)\xi\end{pmatrix},\\
C\begin{pmatrix}f\\g\end{pmatrix}&=f(0),\\
D\xi&=\Theta(0)\xi,
\end{split}
\]
with $f\in\mathcal D(\Theta)$ and $\xi\in\mathbb C^p$, define a
closely connected realization of $\Theta$ which is unitary. This
realization is unique up to a continuous and continuously
invertible similarity operator. \label{tm:unitary}
\end{Tm}

It is important to note that, in some cases, all three
realizations are unitary. This is in particular the case when
$\Theta$ is rational and $J$-unitary on the unit circle. See
Section \ref{fddbs}.

\subsection{Finite dimensional de Branges spaces}
\label{fddbs}
The finite dimensional case is of special importance, and the
case $J=I$ was considered in details in our previous work
\cite{ajl1}. Then the three realizations are unitary, and it is
easier to focus on the $\mathcal P(\Theta)$ spaces. As proved in
\cite[Corollary p. 111]{ad3} for the case $J=I$ and in
\cite[Theorem 5.5, p. 112]{ad3} for the general case, given
$\Theta\in\mathscr S_\kappa^J$, the associated space $\mathcal
P(\Theta)$ is finite dimensional if and only if $\Theta$ is
rational and $J$ unitary on the unit circle:
\[
\Theta(e^{it})^*J\Theta(e^{it})=J,
\]
at all points $e^{it}$ ($t\in [0,2\pi]$) where it is defined. If
moreover $\Theta$ is analytic in a neighborhood of the closed
unit disk, we have
\[
\mathcal P(\Theta)=\mathbf H_{2,J}\ominus\Theta\mathbf H_{2,J}.
\]
Rationality is not enough to insure that $\mathcal P(\Theta)$ is
finite dimensional, as illustrated by the case $J=1$ and
$\Theta=0$. Then, $\mathcal P(\Theta)=\mathbf H_2(\mathbb D)$.\\

\begin{Dn}
We will denote by $\mathscr U_\kappa^J$ the multiplicative group
of rational $\mathbb C^{p\times p}$-valued functions $\Theta$
which take $J$-unitary values on the unit circle, and for which
the corresponding kernel $K_\Theta$ has $\kappa$ negative
squares. We set
\[
\mathscr U^J=\bigcup_{\kappa=0}^\infty \mathscr U_\kappa^J.
\]
\end{Dn}

The results and realizations presented in the previous section
take now an easier form. The various operators can be seen as
matrices. Unitarity above is with respect to the indefinite
metric of $\mathcal P(\Theta)\oplus\mathbb C_J$, and we can
rephrase Theorem \ref{tm:unitary} as follows:

\begin{Tm} Let $W$ be  a rational $\mathbb C^{p\times p}$-valued
function analytic at the origin, and let
\[
W(z)=D+zC(I-zA)^{-1}B
\]
be a minimal realization of $W$. Then, $W$ is $J$-unitary on the
unit circle if and only if there exists an invertible Hermitian
matrix $H$ (which is uniquely determined from the given
realization) such that
\begin{equation}
\label{Stein}
\begin{pmatrix}
A&B\\C&D\end{pmatrix}^*\begin{pmatrix}
H&0\\0&J\end{pmatrix}\begin{pmatrix} A&B\\C&D\end{pmatrix}=
\begin{pmatrix}
H&0\\0&J\end{pmatrix},
\end{equation}
\end{Tm}

The change of variable $z\mapsto 1/z$ yields:

\begin{Tm} Let $W$ be analytic at infinity, and
let
\[
W(z)=D+C(zI-A)^{-1}B.
\]
be a minimal realization of $W$. Then, $W$ is $J$-unitary on the
unit circle if and only if there exists an invertible Hermitian
matrix $H$ (which is uniquely determined from the given
realization) and such that \eqref{Stein} holds.
\end{Tm}

The matrix $H$ is called the {\sl associated Hermitian matrix}
(to the given minimal realization). This result was proved in
\cite[Theorem 3.10]{ag} for the case where $A$ is non-singular.
For the approach using reproducing kernel Hilbert spaces, see
\cite{adrs, ad3, Alpay96}.\\

\section{Cuntz relations}
\label{sec:5}
\subsection{Cuntz relations and the de Branges-Rovnyak spaces}
\label{sec:cuntz}
The results of this section are related to \cite{CMS} and
\cite{ajlm1}. In that last paper, the functions $1,\ldots,
z^{N-1}$ below are replaced by the span of a finite dimensional
backward-shift invariant subspace, but the discussion is
restricted to the Hilbert space case and scalar-valued
functions.\\

Normally by Cuntz relations we refer to a finite system of
isometries $S_1,\ldots, S_N$ in a Hilbert space $\mathcal H$
satisfying two
conditions:\\
$(a)$ Different isometries in the system must have orthogonal
ranges,
\begin{equation}
\label{C1} S_j^*S_k=0,\quad j\not =k,
\end{equation}
and\\
$(b)$ The sum of the ranges equals $\mathcal H$:
\begin{equation}
\label{C2} \sum_{j=1}^N S_jS_j^*=I_{\mathcal H}.
\end{equation}
Note that $(a)$ already forces $\mathcal H$ to be infinite
dimensional. Indeed, if $\mathcal H$ is finite dimensional, an
isometry is unitary and the orthogonality of the ranges is not
possible, see the discussion below and Section
\ref{cuntzrational}. If we allow the isometries to operate
between two finite dimensional spaces of different dimensions,
then one can find isometries which satisfy the Cuntz relations.
It is the set of three conditions: Each $S_i$ is isometric in a
Hilbert space $\mathcal H$, and $(a)$ and $(b)$, together imply
that every realization is a representation of a simple, purely
infinite $C^*$-algebra, called $O_N$.
In applications to filters, the $N$ individual subspaces
represent frequency bands. This allows for versatile
computational algorithms tailored to multiscale problems such as
wavelet decompositions, and analysis on fractals. In our present
paper, we relax some of the original very restrictive axioms,
while maintaining the computational favorable properties. Our
more general framework still allows for algorithms based on
iteration of the operator
family $S_1,\ldots, S_N$ in a particular representation.\\

If one allows isometries between two Hilbert spaces, then the
finite dimensional case may occur, as illustrated by the following
example:
\[
\mathcal H_1=\mathbb C,\quad \mathcal H_2=\mathbb C^2,
\]
and
\[
S_1=\begin{pmatrix}1\\ 0\end{pmatrix},\quad S_2=
\begin{pmatrix}0\\ 1\end{pmatrix}.
\]
We have
\[
S_1^*S_1=S_2^*S_2=1,\quad S_1^*S_2=S_2^*S_1=0,
\]
and
\[
S_1S_1^*+S_2S_2^*=\begin{pmatrix} 1&0\\0& 1\end{pmatrix}.
\]
We go beyond the setting of Hilbert space, and relax the
conditions $(a)$ and $(b)$ imposed in the original framework from
$C^*$-algebra theory, allowing here isometric operators between
two Pontryagin spaces. We still preserve the features of the
representations of use in iterative algorithms.

It is not surprising that in Section \ref{cuntzrational} we have
finite dimensional spaces. Now for the generalized theory, we
must allow for de Branges and Rovnyak spaces, and for negative
squares and signature matrix. The resulting modifications in the
form of the Cuntz relations, in the case of Hilbert space,
entails some non-trivial modifications addressed in the next two
sections. Our main results for this are proved in the present
section, and in Section \ref{cuntzrational}
for the finite dimensional case.\\

The main result of this section is that one can associate in a
natural way to an element $\Theta\in\mathscr S_{k}^J(\mathbb D)$
a family of operators which satisfy the Cuntz relations. We begin
with a preliminary result, which is a corollary of Theorem
\ref{tm:varphi} with $\varphi(z)=z^N$.

\begin{Pn}
Let $\Theta\in \mathscr S_\kappa^J(\mathbb D)$, and let $\mathcal
P(\Theta)$ be the associated Pontryagin space, with reproducing
kernel
\[
K_\Theta(z,w)=\frac{J-\Theta(z)J\Theta(w)^*}{1-zw^*}.
\]
The function
\[
K_\Theta(z^N,w^N)=\frac{J-\Theta(z^N)J\Theta(w^N)^*}{1-z^Nw^{*N}}
\]
has also $\kappa$ negative squares in its domain of definition in
$\mathbb D$. The associated reproducing kernel Pontryagin space $\mathcal M_N$
is equal to the space of functions of the form $F(z)=f(z^N)$,
where $f\in\mathcal P(\Theta)$, with the following indefinite
inner product
\begin{equation}
\label{hermitian}
[ F,G]_{\mathcal M_N}=[ f,g]_{\mathcal P(\Theta)},
\end{equation}
where $g\in \mathcal P(\Theta)$ and $G(z)=g(z^N)$.
\end{Pn}

We have:

\begin{Tm}
Let $\Theta\in \mathscr S_\kappa^J(\mathbb D)$, and let $\mathcal
P(\Theta)$ be the associated Pontryagin space with reproducing
kernel
\[
K_\Theta(z,w)=\frac{J-\Theta(z)J\Theta(w)^*}{1-zw^*}.
\]
Then, for $N\in\mathbb N$, the function $\Theta_N$ defined by
$\Theta_N(z)=\Theta(z^N)$ belongs to $\mathscr S_{N\kappa}^J$.
Furthermore, $\mathcal P(\Theta_N)$ consists of all the functions
of the form
\[
f(z)=\sum_{j=0}^{N-1}z^j f_j(z^N),\quad f_j\in\mathcal P(\Theta).
\]
Any such representation is unique, and the inner product in
$\mathcal P(\Theta_N)$ is given by
\[
[ f,g]_{\mathcal P(\Theta_N)}=\sum_{j=0}^{N-1}[ f_j,g_j]_{\mathcal
P(\Theta)},
\]
where $g(z)=\sum_{j=0}^{N-1}z^jg_j(z^N)$ for some $g_0,\ldots,
g_{N-1}\in\mathcal P(\Theta)$. \label{tm:thetaN}
\end{Tm}

{\bf Proof:} We proceed in a number of steps.\\

STEP 1: {\sl It holds that $\nu_-(\Theta_N)\le N\cdot\kappa$.}\\

Indeed,
\[
\begin{split}
\frac{J-\Theta(z^N)J\Theta(w^N)^*}{1-zw^*}&=
\frac{J-\Theta(z^N)J\Theta(w^N)^*}{1-z^Nw^{N*}}\cdot\frac{1-z^Nw^{N*}}{1-zw^{*}}\\
&=\frac{J-\Theta(z^N)J\Theta(w^N)^*}{1-z^Nw^{N*}}\cdot
(\sum_{k=0}^{N-1}z^kw^{*k}).
\end{split}
\]
This expresses the kernel $K_{\Theta_N}$ as the sum of $N$
kernels, each with $\kappa$ negative squares. Thus,
$\nu_-(\Theta_N)\le N\kappa$. To show that there is equality, we
need to show that the associated spaces have pairwise
intersections which all reduce to the zero function.\\

STEP 2: {\sl Let $k,\ell\in\left\{0,\ldots, N-1\right\}$, such
that $k\not=\ell$. Then, with $\mathcal M_N$ as in the previous theorem:
\[
z^k\mathcal M_N\cap z^\ell \mathcal M_N=\left\{0\right\}.
\]
}
Indeed, assume that $k>\ell$ and let $f,g\in\mathcal M_N$ be such
that
\[
z^kf(z^N)=z^\ell g(z^N).
\]
Then, $f$ and $g$ will simultaneously be identically equal to
$0_{p\times 1}$. Assume $f\not\equiv 0_{p\times 1}$. One of its
components, say the first, with
\mbox{$f=\begin{pmatrix}x_1(z)&\cdots& x_p(z)\end{pmatrix}^t$} is
not identically equal to zero ($p$ is the size of the signature
matrix $J$). Then we obtain
\[
z^{k-\ell}=\frac{y_1(z^N)}{x_1(z^N)},
\]
where $y_1$ denotes the first component of $g$. Since
  $f$ and $g$ are meromorphic in $\mathbb D$, the function $y_1/x_1$ has
  a Laurent expansion at the origin. Moreover the Laurent expansion
  of $\frac{y_1}{x_1}(z^N)$ contains only powers which are multiple of $N$.
By the uniqueness of the Laurent expansion, this contradicts the
fact that it is equal to $z^{k-\ell}$, with $|k-\ell|<N$.\\

STEP 3: {\sl It holds that
\[
\mathcal P(\Theta_N)=\oplus_{j=0}^{N-1}z^j\mathcal M_N,
\]
and it holds that $\nu_{\Theta_N}=N\kappa$.}\\

This is because the spaces $z^j\mathcal M_N$ have pairwise
intersections which reduce to the zero functions in view of STEP
2.
\mbox{}\qed\mbox{}\\

\begin{Tm}
In the notation above, set
\[
(S_jf)(z)=z^jf(z^N)\quad \mathcal
P(\Theta)\longrightarrow\mathcal P(\Theta_N).
\]
Then,
\begin{equation}
S_j^{[*]}f=f_j\quad \mathcal P(\Theta_N)\longrightarrow\mathcal
P(\Theta), \label{sj*}
\end{equation}
and
\begin{equation}
\begin{split}
S_j^{[*]}S_k&=\delta_{j,k}I_{\mathcal P(\Theta)}\\
\sum_{j=0}^{N-1}S_jS_j^{[*]}&=I_{\mathcal P(\Theta_N)},
\end{split}
\label{eq:cuntz}
\end{equation}
where the $[*]$ denotes adjoint between Pontryagin spaces.
\label{th:cuntz12}
\end{Tm}
{\bf Proof:}  We proceed in a number of steps.\\

STEP 1: {\sl The operators $S_j$ are continuous.}\\

The operators $S_j$ are between Pontryagin spaces of different
indices, and some care is required to check continuity. To this
end, fix $j\in\left\{0,\ldots, N-1\right\}$ and note that $S_j$
is everywhere defined. Furthermore we claim that it is a closed
operator. Indeed, let $f_1,f_2\ldots$ be a sequence of elements
in $\mathcal P(\Theta)$ converging strongly to $f\in\mathcal P(\Theta)$
and such that the sequence $S_jf_1, S_jf_2,\ldots$ converges strongly to
$g\in\mathcal P(\Theta_N)$. Strong convergence in a Pontryagin
space implies weak convergence, and in a reproducing kernel
Pontryagin space, weak convergence implies pointwise convergence.
Therefore, for every $w$ where $\Theta$ is defined,
\[
\lim_{k\rightarrow\infty} f_k(w)=f(w),
\]
and
\[
\lim_{k\rightarrow\infty} (S_jf_k)(w)=g(w).
\]
Since $(S_jf_k)(w)=w^jf_k(w)$, and thus $g(w)=w^jf(w)$. Therefore
$g=S_jf$, and the operator $S_j$ is closed, and hence
continuous.\\

STEP 2: {\sl \eqref{sj*} is in force.}\\

Let $g(z)=\sum_{k=0}^{N-1}z^kg_k(z^N)\in\mathcal P(\Theta_N)$
 where the $g_k\in\mathcal P(\Theta)$, and let $u\in
\mathcal P(\Theta)$. Then,
\[
\begin{split}
[S_ju\,,\, g]_{\mathcal P(\Theta_N)}&=[z^ju(z^N)\,,\,
\sum_{k=0}^{N-1}z^kg_k(z^N)]_{\mathcal P(\Theta_N)}\\
&=[u\, ,\, g_j]_{\mathcal P(\Theta)}\\
&=[u\,,\, S_j^{[*]}g]_{\mathcal P(\Theta)},
\end{split}
\]
where $[\,,\,]_{\mathcal P(\Theta)}$ and $[\,,\,]_{\mathcal
P(\Theta_N)}$ denote the indefinite inner products in the
corresponding spaces. Hence, we have $S_j^{[*]}g=g_j$.\\

STEP 3: {\sl The Cuntz relations hold.}\\

From \eqref{sj*} we have for $u\in\mathcal P(\Theta)$
\[
S_j^{[*]}S_ku=S_j^{[*]}(z^ku(z^N))=\begin{cases}0\quad{\rm
if}\quad j\not =k,\\
u\quad{\rm if}\quad j=k.
\end{cases}
\]
Furthermore, for $f(z)=\sum_{j=0}^{N-1}z^jf_j(z^N)\in\mathcal
P(\Theta_N)$ (where the \mbox{$f_j\in\mathcal P(\Theta)$}), we
have
\[
S_kS_k^{[*]}f=S_k(f_k)=z^kf_k(z^N),
\]
and thus
\[
\sum_{k=0}^{N-1}S_kS_k^{[*]}=I_{\mathcal P(\Theta_N)}.
\]
\mbox{}\qed\mbox{}\\

We note that, with
\[
S=\begin{pmatrix}S_0&S_1&\cdots &S_{N-1}\end{pmatrix}\quad
\mathcal P(\Theta)^N\longrightarrow\mathcal P(\Theta_N),
\]
the Cuntz relations \eqref{sj*} can be rewritten as
\[
SS^{[*]}=I_{\mathcal P(\Theta_N)}\quad{\rm and}\quad S^{[*]}S=
I_{\mathcal P(\Theta)^N}.
\]

At this stage, let us introduce some more notation. We set
\[
\Theta_{N^k}(z)=\Theta(z^{N^k}),
\]
and $S_i^{(0)}=S_i$ for $i=0,\ldots, N-1$. We can reiterate the
preceding analysis with $\Theta_N$ instead of $\Theta$. We then
obtain $N$ isometries $S_0^{(1)},\ldots, S_{N-1}^{(1)}$ from
$\mathcal P(\Theta_N)$ into $\mathcal P(\Theta_{N^2})$ satisfying
the Cuntz relations. Iterating $k$ times, one obtains $k$ sets of
isometries,
\[
S_0^{(j-1)},\ldots,S_{N-1}^{(j-1)},\quad j=1,\ldots, k,
\]
from $\mathcal P(\Theta_{N^{j-1}})$ into $\mathcal
P(\Theta_{N^j})$, which also satisfy the Cuntz relations. This
gives us $N^k$ isometries
\[
S_{i_1}^{(0)}S_{i_2}^{(1)}\cdots S_{i_k}^{(k-1)},
\]
with $(i_1,i_2,\ldots ,i_{k})\in\left\{0,\ldots, N-1\right\}^k$,
from $\mathcal P(\Theta)$ into $\mathcal P(\Theta_{N^k})$, all
satisfying the Cuntz relations.
\subsection{Cuntz relation: The general case}
We now wish to extend the results of Section \ref{sec:cuntz}, and
in particular Theorem \ref{tm:thetaN} to the case where the $N$
functions $1,z,\ldots, z^{N-1}$ are replaced by prescribed
functions $m_0(z),m_1(z),\ldots, m_{N-1}(z)$, whose finite
dimensional linear span we denote by $\mathcal L$, and the kernel $K_\Theta(z,w)$ is replaced 
by a given analytic $\mathbb C^{N\times N}$-valued kernel $K(z,w)$ and the kernel
$K_{\Theta_N}(z,w)$ is replaced by a kernel $\widetilde{K}(z,w)$. Let as
in Section \ref{sec:cuntz}, $K_N(z,w)=K(z^N,w^N)$. We 
address the following problem: Given $K$ and $\widetilde{K}$ two
Hermitian kernels defined on a set $\Omega$, and with a finite
number of negative squares there, when can one find decompositions of the form
\begin{equation}
\label{gleason} 
f(z)=\sum_{n=0}^{N-1}m_n(z)g_n(z^N).
\end{equation}
where  the functions $g_0,\ldots, g_{N-1}$ belong to $\mathcal P(K)$
for some, or all, elements in $\mathcal P(\widetilde{K})$.
We have:

\begin{Tm}
Let $K(z,w)$ and $\widetilde{K}(z,w)$ be two kernels defined on a set
$\Omega$, and assume that
\begin{equation}
\label{erty} \nu_-(\widetilde{K})=N\nu_-(K).
\end{equation}
Let $m_0,\ldots, m_{N-1}$ be $N$ functions on $\Omega$. Assume
that the kernel
\[
\widetilde{K}(z,w)-(\sum_{n=0}^{N-1}m_n(z)m_n(w)^*)K(z,w)
\]
is positive definite in $\Omega$. Then, with $\varphi(z)=z^N$,
the choice \mbox{$g_n=T_{m_n,\varphi}^{[*]}f_n$}, $n=0,\ldots,
N-1$ solves \eqref{gleason}. \label{tm:gleason}
\end{Tm}

{\bf Proof:} We use Theorem \ref{th:contraction} with
$K_2(z,w)=\widetilde{K}(z,w)$ and
\[
\quad K_1(z,w)=\begin{pmatrix}
K(z,w)&0&0&\cdots&0\\
0&K(z,w)&0&\cdots&0\\
&& & &\\
& & & &\\
0&0&\cdots&0&K(z,w)\end{pmatrix}.
\]
Then
\[
(\sum_{n=0}^{N-1}m_n(z)m_n(w)^*)K(z,w)=m(z)K_1(z)m(w)^*,
\]
and Theorem \ref{th:contraction} with $K_1$ and $K_2$ as above,
and
\[
m(z)=\begin{pmatrix}m_0(z)&m_1(z)&\cdots&m_{N-1}(z)\end{pmatrix},\quad{\rm
and}\quad\varphi(z)=z^N,
\]
leads to the fact that the map
\[
f\mapsto m(z)f(z^N)
\]
is a contraction from $(\mathcal P(K))^N$ into $\mathcal
P(\widetilde{K})$.
\mbox{}\qed\mbox{}\\

In applications, one uses the kernel $\widetilde{K}(z,w)=K_N(z,w)$ in the above result.

\begin{Pn}
A sufficient condition for \eqref{erty} to hold is that
\begin{equation}
m_j{\mathcal P}(K)\cap m_k\mathcal P(K)=\left\{0\right\},
\end{equation}
for all $j,k\in \left\{0,\ldots, N-1\right\}$ such that $j\not
=k$.
\end{Pn}

{\bf Proof:} Indeed, when this condition is in force, we have that
the Pontryagin space with reproducing kernel $m(z)K_1(z,w)m(w)^*$
is the direct sum of the Pontryagin spaces with reproducing kernels
$m_j(z)K(z,w)m_j(w)^*$, $j=0,\ldots, N-1$.
\mbox{}\qed\mbox{}\\

We note that there are similarity between \eqref{gleason} and the
solution of Gleason's problem: Gleason's problem is the
following: Given a linear space of functions $\mathcal M$ of
functions analytic in a set \mbox{$\Omega\subset\mathbb C^N$},
and given $a\in\Omega$, Gleason's problem is the following: when
can we find functions $g_1,\ldots, g_N\in\mathcal M$ (which
depend on $a$) such that
\[
f(z)-f(a)=\sum_{n=1}^N(z_n-a_n)g_n(z,a)
\]

\subsection{Cuntz relations: Realizations in the rational case}
\label{cuntzrational} Recall that for a
given generalized Schur function $\Theta$, we presented in
Theorems \ref{tm:coiso}  and \ref{tm:unitary}
coisometric and unitary realizations respectively. The unitary
realization turns to be more involved than the coisometric
backwards shift realization. In some cases, these two
realizations are unitarily equivalent, in particular when
$\Theta$ is rational and $J$-unitary on the unit circle. As we
already discussed in Section \ref{fddbs}, this is equivalent to
having the space $\mathcal P(\Theta)$ finite dimensional. In this
section we adopt this simplifying assumption and study the
realization of $\Theta_N(z)=\Theta(z^N)$ in terms of the
realization of $\Theta$.\\

We take the signature matrix $J$ to belong to $\mathbb C^{L\times
L}$. We know (see \cite{ad3, adrs} and Theorem \ref{tm:coiso}
above) that
\[
\Theta(z^N)=\mathscr D+z\mathscr C(I-z\mathscr A)^{-1}\mathscr B
\]
where $\mathscr D=\Theta_N(0)=\Theta(0)$ and  where $\mathscr A,
\mathscr B$ and $\mathscr C$ are defined as follows: $\mathscr C$
is the evaluation at the origin,
\[
\mathscr Cf=f(0).
\]
$\mathscr B$ is defined by
\[
\mathscr B\xi=\frac{\Theta_N(z)-\Theta_N(0)}{z}\xi,\quad
\xi\in\mathbb C^L,
\]
and $\mathscr A$ is the backward shift in $\mathcal P(\Theta_N)$.
The matrix, see \cite{ad3},
\[
\begin{pmatrix}\mathscr A&\mathscr B\\ \mathscr C&\mathscr
D\end{pmatrix}
\]
is unitary in the $\mathcal P(\Theta_N)$ metric. We know from
Theorem \ref{tm:thetaN} that $\mathcal P(\Theta_N)$ is equal to
the space of functions of the form
\begin{equation}
\label{charonne} f(z)=\sum_{k=0}^{N-1}z^kf_k(z^N),
\end{equation}
where the $f_k\in\mathcal P(\Theta)$ are uniquely defined. We
will denote by $U$ the map
\[
f\hookrightarrow \begin{pmatrix} f_0\\ f_1\\ \vdots\\
f_{N-1}\end{pmatrix}
\]
from $\mathcal P(\Theta_N)$ onto $(\mathcal P(\Theta))^N$. In view
of \eqref{eq:cuntz}, $U$ is a unitary map (between Pontryagin
spaces).\\

Let $T$ denote the following map from $(\mathcal P(\Theta))^N$
into itself defined by
\[
  TUf=\begin{pmatrix}0&I&0&\cdots&0\\
0&0&I&\cdots&0\\
& & & &\\
  & & & &I\\
  R_0&0&0&\cdots &0
\end{pmatrix}\begin{pmatrix}f_0\\ f_1\\ \vdots\\ \vdots \\
f_{N-1}\end{pmatrix}=\begin{pmatrix} f_1\\ f_2\\ \vdots\\ \vdots \\
R_0f_0\end{pmatrix}.
\]
\begin{Pn}
Let $f\in\mathcal P(\Theta_N)$, with representation
\eqref{charonne}. It holds that
\begin{equation}
\label{R_0A} U\mathscr A f=(T Uf)(z^N),
\end{equation}
and it holds that
\begin{equation}
\label{A_inner} \langle \mathscr Af,\mathscr Ag\rangle_{\mathcal
P(\Theta_N)}=\langle TUf, TUg\rangle_{(\mathcal P(\Theta))^N}.
\end{equation}
\end{Pn}

{\bf Proof:} Indeed, with $f$ is of the form \eqref{charonne}, we
have
\[
\mathscr Af(z)=
R_0f(z)=\frac{f(z)-f(0)}{z}=\sum_{k=1}^{N-1}z^{k-1}f_k(z^N)+z^{N-1}
\frac{f_0(z^N)-f_0(0)}{z},
\]
so that $U\mathscr AU^*f$ is equal to
\[
\begin{pmatrix}
f_0\\ f_1\\ \vdots \\ f_{N-1}\end{pmatrix}\mapsto\begin{pmatrix}
f_1\\ f_2\\ \vdots \\ R_0f_0\end{pmatrix},
\]
that is \eqref{R_0A} in in force. Finally \eqref{A_inner} follows
from
the formula for the inner product in $\mathcal P(\Theta_N)$. \mbox{}\qed\mbox{}\\

\begin{Pn}
Let $f\in\mathcal P(\Theta_N)$ with representation
\eqref{charonne}. Then,
\begin{equation}
\label{eq:C} \mathscr C f=C\begin{pmatrix}I_L&0&\cdots
&0\end{pmatrix}Uf,
\end{equation}
where $C$ is the evaluation at the origin in $\mathcal P(\Theta)$.
\end{Pn}

{\bf Proof:} This is clear from
\[
\mathscr Cf=f_0(0)=\begin{pmatrix}C&0&0& \cdots&0\end{pmatrix}\begin{pmatrix}f_0\\
                                            f_1\\ \vdots\\ \vdots \\
                                            f_{N-1}\end{pmatrix}
\]
\mbox{}\qed\mbox{}\\

\begin{Pn}
We have
\[
\mathscr B\xi =z^{N-1}(B\xi)(z^N)
\]
where $B$ is the operator from $\mathbb C^L$ into $\mathcal
P(\Theta)$:
\[
B\xi=R_0\Theta\xi
\]
and we have
\begin{equation}
\langle \mathscr B\xi,\mathscr B\eta\rangle_{\mathcal
P(\Theta_N)} =\langle B\xi,B\eta\rangle_{\mathcal
P(\Theta)},\quad\eta,\xi\in\mathbb C^L. \label{B_inner}
\end{equation}
\end{Pn}
{\bf Proof:} We have
\[
\mathscr
B\xi(z)=R_0\Theta_N\xi(z)=\frac{\Theta(z^N)-\Theta(0)}{z}=z^{N-1}(B\xi)(z^N)
\]
Equality \eqref{B_inner} follows form the definition of the inner
product in $\mathcal P(\Theta_N)$.
\mbox{}\qed\mbox{}\\

These various formulas allow to show directly that the
realization is indeed unitary, and to compute the associated
Hermitian matrix in the finite dimensional case.

\section{Decompositions}
\label{sec:6}
\subsection{Generalized down-sampling and an Hermitian form}
In the preceding section we considered decompositions of a
function in the form \eqref{gleason}. Here we consider different
kind of decompositions. We consider matrices $P\in\mathbb C^{N\times N}$ 
satisfying
\begin{equation}
\label{Jussieu_ligne_10}
\det(I_N-\e^\ell P^\ell)\not=0,\quad
\ell=1,\ldots, N-1,\quad{\rm and}\quad P^N=I_N.
\end{equation}
We do not assume that $P^{N-1}\not =I_N$, and in particular the
choice $P=I_N$ is allowed. The special case  $P=\e
P_N$  plays also an important role.

\begin{Tm}
\label{tm:sum} Let $W$ be a $\mathbb C^{N\times M}$-valued
function defined in the open unit disk (typically, $M=1$ or
$M=N$). Let $P\in\mathbb C^{N\times N}$ satisfying
\eqref{Jussieu_ligne_10}, and let, for $k=0,\ldots, N-1$,
\begin{equation}
W_k(z)=\frac{1}{N}\sum_{\ell=0}^{N-1} (\e P)^{k\ell}W(\e^\ell z).
\end{equation}
Then,
\begin{eqnarray}
\label{sym1}
W_k(\e z)&=&(\e P)^{-k}(W_k(z)),\quad k=0,\ldots, N-1,\\
W(z)&=&\sum_{k=0}^{N-1}W_k(z). \label{sum_W}
\end{eqnarray}
\end{Tm}

{\bf Proof:} We have
\[
\begin{split}
W_k(\e z)&=\frac{1}{N}\sum_{\ell=0}^{N-1} (\e P)^{k\ell}W(\e^\ell \e z)\\
&=(\e P)^{-k}\left(\frac{1}{N}\sum_{\ell=0}^{N-1} (\e
P)^{k(\ell+1)}W(\e^{\ell+1} z)\right)
\\
&=(\e P)^{-k}(W_k(z)),
\end{split}
\]
since $(\e P)^{kN}=I_N$, and this proves \eqref{sym1}. To prove 
\eqref{sum_W} we write
\[
\begin{split}
\sum_{k=0}^{N-1}W_k(z)&=\sum_{k=0}^{N-1}\left(\frac{1}{N}\sum_{\ell=0}^{N-1}
(\e P)^{k\ell}W(\e^\ell z)
\right)\\
&=\frac{1}{N}\left(\sum_{\ell=0}^{N-1}
\left(\sum_{k=0}^{N-1}(\e P)^{k\ell}\right)W(\e^\ell z)\right)\\
&=W(z),
\end{split}
\]
since, in view of \eqref{Jussieu_ligne_10},
\[
\sum_{k=0}^{N-1}(\e P)^{k\ell}=\begin{cases}
N,\,\,{\rm if}\, \ell=0,\\
(I_N-(\e P)^{N\ell})(I-(\e P)^{\ell})^{-1}=0\, {\rm if}\,
\ell=1,2,\ldots N-1.\end{cases}
\]
\mbox{}\qed\mbox{}\\

When $P=I_N$, the index $k=1$ corresponds to the down-sampling
operator.

\subsection{Orthogonal decompositions in Krein spaces}
\label{azxcv}
In some cases the decomposition \eqref{sum_W} is orthogonal for
the underlying Krein space (or Pontryagin space) structure. We
will assume that the Krein space $(\mathcal K,
[\cdot,\cdot]_{\mathcal K})$ consists of $\mathbb C^N$-valued
functions and satisfies the following property:
\begin{Hyp} Let $P$ be a matrix satisfying
\eqref{Jussieu_ligne_10}, and let $\varphi(z)=\e z$. We assume
that:\\
$(1)$ The composition operator $f\mapsto f(\varphi)$ is continuous and unitary 
from $\mathcal K$ into itself.\\ 
$(2)$ The operator of multiplication by $P$ on the left 
is continuous and unitary from $\mathcal K$ into itself.
\label{hyp}
\end{Hyp}
We note that, in particular, the operator $T_{P,\varphi}$ defined by  \eqref{tm},
\[
T_{P,\varphi}f(z)=Pf(\e z),
\]
is continuous and unitary from $\mathcal K$ into itself. Note also that
\[
T_{P,\varphi}^N=I_{\mathcal K}.
\]

Hypothesis \ref{hyp} hold in particular for the spaces $\mathbf
H_{2,J}$ when $P$ is $J$-unitary, that is, satisfies
\[
P^*JP=J.
\]
\begin{Tm}
Let $(\mathcal K, [\cdot\, ,\,\cdot\, ])$ be a Krein space of
$\mathbb C^{N}$-valued functions, satisfying Hypothesis
\ref{hyp}. Let $W\in\mathcal K$ and let
\begin{equation}
W_k(z)=\frac{1}{N}\sum_{\ell=0}^{N-1}(\e P)^{k\ell}W(\e^\ell z).
\end{equation}
Then,
\[
[W_\ell,W_k]=0,\quad \ell\not= k,
\]
\[
W(z)=W_0(z)+\cdots+W_{N-1}(z),
\]
and
\[
W_k(\e z)=(\e P)^{-k}W(z).
\]
\end{Tm}

{\bf Proof:} The last two claims are proved in Theorem
\ref{tm:sum}. The first claim takes into account the hypothesis
on $\mathcal K$, and is proved as follows: We take $k_1$ and
$k_2$ in $\left\{0,\ldots, N-1\right\}$, and assume that
$k_2<k_1$. Taking into account the definition of $W_k$, we see
that the inner product $[W_{k_1},W_{k_2}]_{\mathcal K}$ is a sum
of $N^2$ inner products, namely
\[
[(\e P)^{k_1 \ell_1}W(\e^{\ell_1}z),(\e
P)^{k_2\ell_2}W(\e^{\ell_2}z)]_{\mathcal K},\quad
\ell_1,\ell_2\in\left\{0,\ldots, N-1\right\}.
\]
These $N^2$ inner products can be rearranged as $N$ sums of inner
product, each sum being equal to $0$. Indeed, consider first the
inner products corresponding to $\ell_1=\ell_2$. In view of the
unitary of the operator $T_{P,\varphi}$ we have
\[
\begin{split}
\sum_{\ell_1=0}^{N-1}[(\e P)^{k_1 \ell_1}W(\e^{\ell_1}z),(\e P)
^{k_2\ell_1}W(\e^{\ell_1}z)]_{\mathcal
K}&=[\left(\sum_{\ell_1=0}^{N-1}(\e^{k_1-k_2} P
)^{\ell_1}\right)W,W]_{\mathcal
K}\\
&=0.
\end{split}
\]
Indeed, using $0<k_1-k_2\le N-1$, and so, by hypothesis on $P$, we
have
\[
\det(I_N-(\e P)^{k_1-k_2})\not =1,
\]
and the sum
\[
\sum_{\ell_1=0}^{N-1}((\e P)^{k_1-k_2})^{\ell_1}=0.
\]
Let us now regroup the factors of $[W(z),W(\e z)]_{\mathcal K}$.
Taking into account that 
\[
[P^{k_1(N-1)}W(\e^{N-1}z),W(z)]_{\mathcal K}=[P^{k_1(N-1)}W(z), W(\e z)]_{\mathcal
K},
\]
we have
\[
\begin{split}
\sum_{\ell=0}^{N-2}[(\e P)^{k_1\ell}W(\e^\ell z),(\e P)
^{k_2(\ell+1)}W(\e^{\ell+1}z)]_{\mathcal
K}+&\\
&\hspace{-9.5cm}+[(\e
P)^{k_1(N-1)}W(\e^{N-1}z),W(z)]_{\mathcal K}
\\
&\hspace{-10cm}=[ \left(\sum_{\ell=0}^{N-2} (\e P)^{\ell
k_1-(\ell+1)k_2}+(\e P)^{k_1(N-1)}\right) W(z),W(\e z)]_{\mathcal
K}\\
&\hspace{-10cm}=[(\e P)^{-k_2}\left(\sum_{\ell=0}^{N-1}((\e
P)^{k_1-k_2})^\ell \right)W(z), W(\e z)]_{\mathcal
K}\\
&\hspace{-10cm}=0.
\end{split}
\]
The remaining terms are summed up to $0$ in the same way.
\mbox{}\qed\mbox{}\\

\subsection{Decompositions in reproducing kernel spaces}
We begin with a result in the setting of Schur functions, as
opposed to generalized Schur functions.
\begin{Tm}
Let $W$ be a $\mathbb C^{p\times q}$-valued Schur function and let
$\varphi(z)=\e z$. Then the operator of composition by $\varphi$
is a contraction from $\mathcal H(W)$ into itself if and only if
there exists a $\mathbb C^{q\times q}$-valued Schur function
$X(z)$ such that
\begin{equation}
\label{republique} W(z) =W(\e z)X(z).
\end{equation}
\label{tm:dec_rkhs}
\end{Tm}

{\bf Proof:} By Theorem \ref{th:contraction}, the map $T_\varphi$
is a contraction if and only if the kernel
\[
K_W(z,w)-K_W(\e z, \e w)=\frac{W(\e z)W(\e
w)^*-W(z)W(w)^*}{1-zw^*}
\]
is positive definite in the open unit disk. By Leech's
factorization theorem, see \cite[p. 107]{rr-univ}, the above
kernel is positive definite if and only if there is a Schur
function $X(z)$ such that \eqref{republique} is in force.
\mbox{}\qed\mbox{}\\

As an example, take any Schur function $s$ and build
\begin{equation}
W(z)=\frac{1}{\sqrt{N}}\begin{pmatrix}s(z)&s(\e z)&\cdots &
s(\e^{N-1}z)\end{pmatrix}. \label{azxc}
\end{equation}
Then
\[
W(z)=W(\e z)P_N,
\]
where $P_N$ is defined by \eqref{Cardinal_Lemoine_Ligne_10}.
\section{The family $\mathscr{C}_N$}
\label{sec:7}
An effective approach to generating wavelet bases is the use of
Multiresolution Analysis (MRA), see for example \cite{BJMP05,
BrJo02a, Dau92}. Traditionally one looks for a finite family of
functions in $\mathbf L_2(\mathbb R,dx)$, or  $\mathbf L_2(\mathbb
R^d,dx)$ for some dimension $d$. If $d=1$, one chooses a scale
number, say $N$. If $d>1$, instead one scales with a  $d\times d$
matrix $A$ over the integers. We assume that $A$ is expansive,
i.e., with eigenvalues bigger than $1$ in modulus. If $A$ is
given, let $N$ be the absolute value of its determinant. To create
MRA wavelets we need an initial finite family $\mathcal F$  of $N$
functions in  $\mathbf L_2(\mathbb R)$, or  $\mathbf L_2(\mathbb
R^d)$. One of the functions is called the scaling function ($\phi$
in the discussion below). For the moment, we will set $d=1$, but
the outline below easily generalizes to $d>1$. An MRA wavelet
basis is a basis for  $\mathbf L_2(\mathbb R)$, or  $\mathbf
L_2(\mathbb R^d)$ which is generated from the initial family
$\mathcal F$ and two operations : one operation is scaling by the
number $N$ (or the matrix $A$ if $d>1$), and the other is action
by integer translates of functions. The special property for the
finite family of functions $\mathcal F$ is that if the $N$-scaling
is applied each function $\psi$ in F the result is in the closed
span of the integer translates of the scaling function $\phi$. The
corresponding coefficients are called masking coefficients. The
reason for this is that the scaled functions represent
refinements, and they are computed from masking points in a
refinement. The role of the functions $m_0, m_1, \ldots  ,
m_{N-1}$ are the frequency response functions corresponding to
the system of masking coefficients. From these functions we then
build a matrix valued function $W(z)$ as in \eqref{eq:m1}. The
question we address here is the characterization of the matrix
valued function which arise this way. Now the wavelet filters we
consider here go beyond those studied earlier in that we allow
for wider families of Multiresolution Analyses (MRAs). This
includes more general wavelet families, allowing for example for
wavelet frame bases, see e.g., \cite{BJMP05, Jor08, JoSo09},
multi-scale systems in dynamics, and in analysis of fractals; see
\cite{DuJo06b}.
\subsection{The family $\mathscr{C}_N$: characterization}
The filters we consider are matrix-valued (or operator valued)
functions of a complex variable. In general if a positive integer
$N$ is given, and if a matrix function $W(z)$ is designed to take
values in $\mathbb C^{N\times N}$, then of course, there are
$N^2$ scalar-valued function occurring as matrix entries.
However, in the case of filters arising in applications
involving  $N$ distinct frequency-bands, for example in wavelet
constructions with scale number $N$, then we can take advantage
of an additional symmetry for the given matrix function $W(z)$,
see for example \eqref{eq:sym11} in the Introduction. Here we
point out that this $N$-symmetry condition (or $N$-periodicity)
means that $W(z)$ is then in fact determined by only $N$ scalar
valued functions, see \eqref{eq:m1} below. These functions play
three distinct roles as follows: They are (i) the scalar valued
filter functions, $\hat{s}_i$, for $i = 0, 1,\ldots, N-1$,  in
generalized quadrature-mirror filter systems (the quadrature case
corresponds to $N =2$); they are (ii) scaling filters for
scale-number $N$ with each of the $N$ scalar functions
$\hat{s}_i$  generating an element in a wavelet system of
functions on the real line and corresponding to scale-number $N$;
and (iii) the system of scalar functions $(\hat{s}_i)_{i=0,\ldots
N-1}$ generates an operator family $( S_i)_{i=0,\ldots, N-1}$
constituting a representation of the Cuntz relations; thus
generalizing Theorem 5.3 above. The results presented in this
section are related to \cite{ajl1}.\\

Recall that $\e=e^{\frac{2\pi i}{N}}$. We shall say that a
$\mathbb C^{N\times N}$-valued ($N\geq 2$) function $W$
meromorphic in the open unit disk $\mathbb D$ belongs to
$\mathcal{C}_N$ if it is of the form
\begin{equation}
\label{eq:m1} W(z)=\frac{1}{\sqrt{N}}\begin{pmatrix} \hat{s}_0(z)&
\hat{s}_0(\e z)&\cdots &\hat{s}_0(\e^{N-1}z)\\
\hat{s}_1(z)&\hat{s}_1(\e z)&\cdots &\hat{s}_1(\e^{N-1}z)\\
\vdots& \vdots&  &\vdots \\
\hat{s}_{N-1}(z)&\hat{s}_{N-1}(\e z)&\cdots &
\hat{s}_{N-1}(\e^{N-1}z)
\end{pmatrix},
\end{equation}
where $\hat{s}_0,\ldots, \hat{s}_{N-1}$ are complex-valued
functions meromorphic in $\mathbb D$. Note that such a function,
when analytic at the origin, will never be invertible there. A
special case of this analyticity restriction of course is when
$W(z)$ has polynomial entries. Under the filter-to-wavelet
\footnote{This correspondence: {\sl polynomial filter} to {\sl
compactly supported wavelet} even works if $d>1$.}correspondence
\cite{BrJo02a}, polynomial filters are the compactly supported
wavelets. In the sequel, it will turn out that we shall
concentrate on the opposite cases. Namely, not only $W(z)$ will
not be analytic at the origin, in fact we shall
have $W(z)^{-1}_{|_{z=0}}=0_{N\times N}$.\\

Recall that we have denoted by $P_N$ the permutation matrix,
\[
P_N=\begin{pmatrix} 0_{1\times (N-1)}&1\\
I_{N-1}& 0_{(N-1)\times 1}\end{pmatrix}
\]
(see \eqref{Cardinal_Lemoine_Ligne_10}).
\begin{La}
A $\mathbb C^{N\times N}$-valued function meromorphic in the open
unit disk is of the form \eqref{eq:m1} if and only if it satisfies
\eqref{eq:sym11}:
\begin{equation*}
W(\e z)=W(z)P_N
\end{equation*}
\label{lemma2.1}
\end{La}

{\bf Proof:} Let $W$ be a $\mathbb C^{N\times N}$-valued function
meromorphic in $\mathbb D$, and satisfying \eqref{eq:sym11}, and
let $s_1,\ldots s_N$ denote its columns, i.e.
\begin{equation}
\label{opera_bastille}
W(z)=\begin{pmatrix}s_1(z)&s_2(z)&\ldots&s_N(z)\end{pmatrix}.
\end{equation}
Namely, from \eqref{eq:m1}
\[
s_j(z):=\frac{1}{\sqrt{N}}\begin{pmatrix}
\hat{s}_0(\e^{j-1}z)\\
\hat{s}_1(\e^{j-1}z)\\
\vdots\\
\hat{s}_{N-1}(\e^{j-1}z)
\end{pmatrix},\quad\quad\quad j=1,~\cdots~,~N.
\]
Multiplying $W$ by $P_N$ from the right makes a cyclic shift of
the columns to the left, namely
\[
W(z)P_N=\begin{pmatrix}s_2(z)&s_3(z)&\cdots
&s_{N}(z)&s_1(z)\end{pmatrix}.
\]
Equation \eqref{eq:sym11} then leads to
\[
\begin{split}
\begin{pmatrix}s_1(\e z)&s_2(\e z)&\cdots
&s_{N-1}(\e z)&s_N(\e z) \end{pmatrix} &=\\
&\hspace{-3.cm}=
\begin{pmatrix}s_2(z)&s_3(z)&\cdots
&s_{N}(z)&s_1(z)\end{pmatrix}.
\end{split}
\]
Thus
\[
s_2(z)=s_1(\e z),\quad s_3(z)=s_1(\e^2 z), \ldots,
s_N(z)=s_1(\e^{N-1}z),
\]
and so $W$ is of the asserted form. The converse is clear.
\mbox{}\qed\mbox{}\\

Note that in contrast to Lemma \ref{lemma2.1}, in equation
\eqref{azxc} we did not assume that $W$ is square.\\

When one assumes that the function $W$ in the previous lemma is a
generalized Schur function, the symmetry condition
\eqref{eq:sym11} can be translated into the realization. We
present the result for the closely outer connected coisometric
realization, but similar results hold for the closely inner
connected isometric realization and connected unitary
realizations as well (see Section \ref{sub42} for these notions).
In the statement, recall that the state space $\mathcal P$ will
in general be infinite dimensional and endowed with a Pontryagin
space structure.

\begin{Tm}
Let $W$ be a generalized Schur function, and let
\[
W(z)=D+zC(I-zA)^{-1}B
\]
be a closely-inner coisometric realization of $W$, with state
space $\mathcal P$. Then, $W$ satisfies \eqref{eq:sym11} if and
only if there is a bounded invertible operator $T$ from $\mathcal
H$ into itself such that

\begin{equation}
\label{uniq}
\begin{pmatrix}
\e A&B\\ \e
C&D\end{pmatrix}\begin{pmatrix}T&0\\0&I_N\end{pmatrix}=
\begin{pmatrix}T&0\\0&I_N\end{pmatrix}
\begin{pmatrix}
 A&B\\  C&D\end{pmatrix}
\end{equation}
Furthermore, the operator $T$ satisfies:
\begin{equation}
T^N=I.
\end{equation}
\end{Tm}

{\bf Proof:} The first equation follows from the uniqueness of
the closely connected coisometric realization. Iterating
\eqref{uniq} and taking into account that $\e^N=1$ we get
\[
\begin{pmatrix}
A&B\\  C&D\end{pmatrix}\begin{pmatrix}T^N&0\\0&I_N\end{pmatrix}=
\begin{pmatrix}T^N&0\\0&I_N\end{pmatrix}
\begin{pmatrix}
 A&B\\  C&D\end{pmatrix}.
 \]
By uniqueness of the similarity operator we have $T^N=I$.
\mbox{}\qed\mbox{}\\
\begin{Pn}
Let $W_1$ and $W_2$ be in $\mathcal C_N$. Then the functions
\[
W_1(z)W_2(\overline{z})^*\quad{and}\quad
W_1(z)W_2(1/\overline{z})^*
\]
are meromorphic functions of $z^N$.
\end{Pn}

{\bf Proof:} Let $W(z)=W_1(z)W_2(\overline{z})^*$. Since
$P_NP_N^*=I_N$, we have
\[
\begin{split}
W(\e z)&=W_1(\e z)W_2(\overline{\e z})^*\\
&=W_1(z)P_NP_N^*W_2(\overline{z})^*\\&=W_1(z)W_2(\overline{z})^*\\&=W(z),
\end{split}
\]
that is
\begin{equation}
W(\e z)=W(z). \label{bastille111}
\end{equation}
The function $W_1$ and $W_2$ are meromorphic in the open unit
disk, and so is the function $W$. We denote by $\Lambda$ the set
of poles of $W$ and by $\Lambda_N$ the set of points $w$ in the
open unit disk such that $w^N\in \Lambda$. Let for
$z=re^{i\theta}$ with $r>0$ and $\theta\in(-\pi,\pi]$,
\[
R(z)=W(\sqrt[N]re^{i\frac{\theta}{N}}).
\]
The function $R$ is analytic in $\mathbb
D\setminus\left\{\Lambda_N\cup(-1,0]\right\}$. Thanks to
\eqref{bastille111} it is continuous across the negative axis at
those points in $(-1,0)$ which are not in $\mathbb D\setminus
\Lambda_N$. It follows that $R$ is analytic in $\mathbb D\setminus
\Lambda_N\cup\left\{0\right\}$. Furthermore, $W(z)=R(z^N)$. Any
singular point of $R$ is a pole (otherwise its roots of order $N$
would be essential singularities of $W$), and so $R$ is
meromorphic in $\mathbb D$.
\mbox{}\qed\mbox{}\\

In the rational case, the previous result has an easier and more
precise proof. Indeed consider the Laurent expansion at the
origin of $W$:
\[
W(z)=\sum_{-m_0}^\infty W_kz^k.
\]
It converges in a punctured disk $0<|z|<r$ for some $r>0$.
Equation \eqref{bastille111} implies that
\[
\sum_{-m_0}^\infty W_kz^k=\sum_{-m_0}^\infty W_k\e^kz^k.
\]
By uniqueness of the Laurent expansion we get that
\[
W_k=0,\quad {\rm for}\quad k\not\in N\mathbb Z.
\]
Thus, if $m>0$, we may assume without loss of generality that
\mbox{$m_0=Nn_0$} for some $n_0\in\mathbb N$. The function
\[
W_-(z)=\sum_{k=-m_0}^{-1} W_kz^k
\]
is rational, and so is the function
\[
W_+(z)=\sum_{k=0}^\infty W_kz^k.
\]
We see that
\[
W_-(z)= \sum_{-m_0\le nN\le -N}W_{nN}z^{nN}
\]
and so $W_-(z)=R_-(z^N)$, where the function
\[
R_-(z)=\sum_{-m_0\le nN\le -N} W_{nN}z^n
\]
is rational and analytic at infinity. The function $W_+$ is
analytic at the origin, and thus can be written in realized form
as:
\[
W_+(z)=D+zC(I_p-zA)^{-1}B.
\]
Comparing with
\[
W_+(z)=\sum_{n=0}^\infty W_{nN}z^{nN}
\]
we have that
\[
CA^pB=\begin{cases}0\quad\hspace{.6cm}{\rm if}\quad
p+1\not\in N\mathbb N,\\
          W_{nN}\quad{\rm if}\quad p+1=nN, \quad n\in\mathbb N.
          \end{cases}
          \]It follows that $W_+(z)=R_+(z^N)$, where $R_+$ is the rational
function defined by
\[
\begin{split}
R_+(z)&=D+\sum_{n=1}^\infty z^{n}CA^{nN-1}B\\
&=D+\sum_{n=1}^\infty z^nCA^{(n-1)N}A^{N-1}B\\
&=D+zC(I_p-zA^N)^{-1}A^{N-1}B.
\end{split}
\]
The function
\[
R(z)=R_-(z)+R_+(z)
\]
is rational.\\

The proof of  the preceding proposition can be mimicked to obtain
the following result:

\begin{Pn}
\label{MR} Let $W_1$ and $W_2$ be in $\mathcal C_{N}$, with non
identically vanishing determinant. Then there exists a
meromorphic function $R$ such that
\begin{equation}
\label{eq7.4} W_1(z)W_2(z)^{-1}=R(z^N).
\end{equation}
\end{Pn}

  To this end, recall that the unitary
matrix $F_N$,
\begin{equation*}
F_N:=\frac{1}{\sqrt{N}}
\begin{pmatrix}
\e^{-(0\cdot 0)}&\e^{-(0\cdot 1)}&\e^{-(0\cdot 2)}&\cdots
&e^{-(0\cdot(N-1))}\\
\e^{-(1\cdot 0)}&\e^{-(1\cdot 1)}&\e^{-(1\cdot 2)}&
\cdots&e^{-(1\cdot(N-1)}\\
\e^{-(2\cdot 0)}&\e^{-(2\cdot 1)}&\e^{-(2\cdot 2)}&
\cdots&e^{-(2\cdot(N-1)}\\
\vdots& \vdots& \vdots& \vdots&\vdots\\
\e^{-((N-1)\cdot 0)}&\e^{-((N-1)\cdot 1)}&\e^{-((N-1)\cdot 2)}&
\cdots&e^{-((N-1)\cdot(N-1)}
\end{pmatrix}.
\end{equation*}
generates the discrete Fourier transform. Namely, the discrete
Fourier transform of $x\in\C^N$ is given by $X=F_Nx$ and the
inverse  discrete Fourier transform is given by $x=F_N^*X$. Let
furthermore
\begin{equation}
\label{eq3.6} \hat{W}_N(z):={\rm
diag}\{1,~z^{-1},~\ldots~,~z^{1-N}\}F_N.
\end{equation}

With this special choice of $W_2$ the previous proposition
becomes:
\begin{Pn}\label{FamilyW}
$W\in\mathcal{C}_N$ and $\det W\not\equiv0$ if and only if it can
be written as
\[
W(z)=R(z^N)\hat{W}_N(z),
\]
where $R$ and $\hat{W}_N$ are as in \eqref{eq7.4} and
\eqref{eq3.6} respectively.
\end{Pn}

\subsection{A connection with periodic systems}
Let
\[
D_N(z)={\rm diag}~(z^N,z^{N-1}\e^{N-1},z^{N-2}\e^{k-2},\ldots,
z\e),
\]
so that
\[
D_N(1)={\rm diag}~(1,\e^{N-1},\e^{k-2},\ldots, \e).
\]
Functions which satisfy the related symmetry
\begin{equation}
\label{eq:abl4}
W(\e z)=D_N(1)^{-1}W(z)P_N
\end{equation}
appear in the theory of periodic systems. A function $W$
satisfies \eqref{eq:abl4} if and only if it is of the form
\begin{equation}
\label{eq:m11}
W(z)=\frac{1}{\sqrt{N}}\begin{pmatrix}
\hat{s}_0(z)&
\hat{s}_0(\e z)&\cdots &\hat{s}_0(\e^{N-1}z)\\
\hat{s}_1(z)&\frac{1}{\e}\hat{s}_1(\e z)&\cdots
&\frac{1}{\e^{N-1}}\hat{s}_1(\e^{N-1}z)\\
\vdots& &  & \\
\hat{s}_{N-1}(z)&\frac{1}{\e^{N-1}}\hat{s}_{N-1}(\e z)&\cdots &
\frac{1}{\e^{(N-1)^2}}\hat{s}_{N-1}(\e^{N-1}z).
\end{pmatrix}
\end{equation}
See \cite[Theorem 4.1, p. 381]{abl4}. We note that the
corresponding general bitangential interpolation problem (see
\cite{bgr} for references) was solved in \cite{abl4} for
functions analytic and contractive in the open unit disk (that
is, for Schur functions). Let us denote by ${\mathcal Per}_N$ the
family of functions meromorphic in the open unit disk and which
satisfy \eqref{eq:abl4}.

\begin{Pn}
The map $W\mapsto D_NW$ is one-to-one from ${\mathcal
Per}_N$ onto $\mathcal{C}_N$. If $W$ is analytic and contractive
in the open unit disk so is $D_NW$.
\end{Pn}

{\bf Proof:} We first note that
\begin{equation}
\label{DN1}
D_N(\e z)=D_N(z)D_N(1).
\end{equation}
Let now $W\in {\mathcal Per}_N$. In view of \eqref{DN1} and
\eqref{eq:abl4} we have
\[
\begin{split}
D_N(\e z)W(\e z)&=D_N(z)D_N(1)D_N(1)^{-1}W(z)P_N\\
&=D_N(z)W(z)P_N,
\end{split}
\]
and so $D_NW\in\mathcal C_N$.
\mbox{}\qed\mbox{}\\

{\bf Epilogue:} A reason for the recent success of wavelet
algorithms is a coming together of tools from engineering and
harmonic analysis. While wavelets now enter into a multitude of
applications from analysis and probability, it was the
incorporation of ideas from signal processing that offered new
and easy-to-use algorithms, and hence wavelets are now used in
both discrete problems, as well as in harmonic analysis
decompositions. 
Following this philosophy
we here employed tools from system theory to wavelet problems
and tried to show how ideas from wavelet 
decompositions throw light on factorizations used by engineers.\\

Since workers in wavelet theory often are not familiar with
filterers in general, and FIR filters  (short for Finite Impulse Response) 
in particular, widely used
in the engineering literature, we have taken the opportunity to
include a section for mathematicians about filters. Conversely
(in the other direction), engineers are often not familiar with
wavelet analysis, and we have included a brief exposition of
wavelet facts addressed to engineers . We showed that
there are explicit actions of infinite-dimensional Lie groups
which accounts for all the wavelet filters; as well as for other
classes of filters used in systems theory. Moreover, we described these
groups, and explained how they arise in systems. The corresponding
algorithms, including the discrete wavelet algorithms are used in
a variety of multi-scale problems, as used for example in data
mining. These are the discrete algorithms, and we described  their
counterparts in harmonic analysis in standard  $\mathbf L_2$
Lebesgue spaces, as well as in reproducing kernels Hilbert spaces.
We also outlined the role of Pontryagin spaces in the study of
stability questions.\\

In the engineering literature the study of filters is mostly
confined to FIR filters. Recall that FIR filters correspond to having the spectrum at the origin.
In our previous work \cite{ajl1} we have explained that the set of
FIR wavelet filters is small in a sense we made precise. 
This suggests two possible conclusions,

(i) It is unrealistic to offer optimization schemes, over
all FIR wavelet filters as part of the design procedure.

(ii) It calls upon using, at least in some circumstances, also 
stable IIR (short for infinite impulse response)
wavelet filters, i.e. the spectrum is confined to the open
unit disk.\\

The above extension to $\mathcal U^{I_N}$ allows us to consider
filters whose spectrum is in $\mathbb C\setminus\mathbb T$. The 
generalization to $\mathcal U^J$ permits the spectrum to be everywhere 
in the complex plane.\\

Roughly, we hope that this message will be useful to
practitioners in their use of these rigorous mathematics tools. We
offer algorithms hopefully improving on those used before.\\

{\bf Acknowledgments:} {D. Alpay thanks the Earl Katz family for
endowing the chair which supported his research. The work was
done in part while the second named author visited Department of
Mathematics, Ben Gurion University of the Negev, supported by a
BGU distinguished visiting scientist program. Support and
hospitality is much appreciated. We acknowledge discussions with
colleagues there, and in the US, Dorin Dutkay, Myung-Sin Song,
and Erin Pearse.}

\bibliographystyle{plain}
%
\def\cprime{$'$} \def\lfhook#1{\setbox0=\hbox{#1}{\ooalign{\hidewidth
  \lower1.5ex\hbox{'}\hidewidth\crcr\unhbox0}}} \def\cprime{$'$}
  \def\cprime{$'$} \def\cprime{$'$} \def\cprime{$'$} \def\cprime{$'$}

\end{document}